\documentclass[11pt]{amsart}
\usepackage{mathrsfs}
\usepackage{amssymb,latexsym}
\usepackage{pstcol,pstricks,color}
\usepackage{chemarrow}

 \numberwithin{equation}{section}

\newtheorem{theorem}{Theorem}[section]
\newtheorem{proposition}[theorem]{Proposition}
\newtheorem{corollary}[theorem]{Corollary}

\newtheorem{remark}[theorem]{Remark}

\newtheorem{lemma}[theorem]{Lemma}

\def\qed{\hfill $\Box$}
\def\pf{\noindent {\it Proof.} }

\title[]{The $\mathbf{uvu}$-avoiding $(a,b,c)$-Generalized Motzkin paths with vertical steps: bijections and statistics enumerations}

\begin{document}
\maketitle

\begin{center}
Yidong Sun\footnote{Corresponding author: Yidong Sun.}, Weichen Wang$^{2}$ and Cheng Sun$^{3}$

School of Science, Dalian Maritime University, 116026 Dalian, P.R. China\\[5pt]

{\it Emails: $^{1}$sydmath@dlmu.edu.cn, $^{2}$weichenw@dlmu.edu.cn, $^{3}$scmath@dlmu.edu.cn}

\end{center}\vskip0.2cm

\subsection*{Abstract} A generalized Motzkin path, called G-Motzkin path for short, of length $n$ is a lattice path from $(0, 0)$ to $(n, 0)$ in the first quadrant of the XOY-plane that consists of up steps $\mathbf{u}=(1, 1)$, down steps $\mathbf{d}=(1, -1)$, horizontal steps $\mathbf{h}=(1, 0)$ and vertical steps $\mathbf{v}=(0, -1)$. An $(a,b,c)$-G-Motzkin path is a weighted G-Motzkin path such that the $\mathbf{u}$-steps, $\mathbf{h}$-steps, $\mathbf{v}$-steps and $\mathbf{d}$-steps are weighted respectively by $1, a, b$ and $c$.
In this paper, we first give bijections between the set of $\mathbf{uvu}$-avoiding $(a,b,b^2)$-G-Motzkin paths of length $n$ and the set of $(a,b)$-Schr\"{o}der paths as well as the set of $(a+b,b)$-Dyck paths of length $2n$, between the set of $\{\mathbf{uvu, uu}\}$-avoiding $(a,b,b^2)$-G-Motzkin paths of length $n$ and the set of $(a+b,ab)$-Motzkin paths of length $n$, between the set of $\{\mathbf{uvu,uu}\}$-avoiding $(a,b,b^2)$-G-Motzkin paths of length $n+1$ beginning with an $\mathbf{h}$-step weighted by $a$ and the set of $(a,b)$-Dyck paths of length $2n+2$. In the last section, we focus on the enumeration of statistics ``number of $\mathbf{z}$-steps" for $\mathbf{z}\in \{\mathbf{u}, \mathbf{h}, \mathbf{v}, \mathbf{d}\}$ and ``number of points" at given level in $\mathbf{uvu}$-avoiding G-Motzkin paths. These counting results are linked with Riordan arrays.

\medskip

{\bf Keywords}: G-Motzkin path, $(a,b)$-Dyck path, $(a,b)$-Motzkin path, $(a,b)$-Schr\"{o}der path, Riordan array.

\noindent {\sc 2020 Mathematics Subject Classification}: Primary 05A15, 05A19; Secondary 05A10.

{\bf \section{ Introduction } }

A {\it generalized Motzkin path}, called {\it G-Motzkin path} for short, of length $n$ is a lattice path from $(0, 0)$ to $(n, 0)$ in the first quadrant of the XOY-plane that consists of up steps $\mathbf{u}=(1, 1)$, down steps $\mathbf{d}=(1, -1)$, horizontal steps $\mathbf{h}=(1, 0)$ and vertical steps $\mathbf{v}=(0, -1)$.

A {\it $\mathbf{ud}$-peak} ({\it $\mathbf{uv}$-peak}) in a G-Motzkin path is an occurrence of $\mathbf{ud}$ ($\mathbf{uv}$). A {\it $\mathbf{du}$-valley} ({\it $\mathbf{vu}$-valley}) in a G-Motzkin path is an occurrence of $\mathbf{du}$ ($\mathbf{vu}$). A point of a G-Motzkin path with ordinate $\ell$ is said to be at {\it level} $\ell$. A step of a G-Motzkin path is said to be at level $\ell$ if the ordinate of its endpoint is $\ell$. By a {\it return step} we mean a $\mathbf{d}$-step or a $\mathbf{v}$-step at level $0$. A G-Motzkin path $\mathbf{P}$ is said to be {\it primitive} if $\mathbf{P}=\mathbf{u}\mathbf{P}'\mathbf{d}$ or $\mathbf{P}=\mathbf{u}\mathbf{P}'\mathbf{v}$ for certain G-Motzkin path $\mathbf{P}'$. A {\it matching step} of a $\mathbf{u}$-step at level $k\geq 1$ in a G-Motzkin path is the leftmost step among all $\mathbf{d}$-steps or  $\mathbf{v}$-steps at level $k-1$ right to the $\mathbf{u}$-step. See Figure 1 for a G-Motzkin path of length $25$ with two $\mathbf{ud}$-peaks, four $\mathbf{uv}$-peaks, two $\mathbf{du}$-valleys, two $\mathbf{vu}$-valleys and two $\mathbf{uvu}$ strings.

\begin{figure}[h] \setlength{\unitlength}{0.5mm}

\begin{center}
\begin{pspicture}(13,2.2)
\psset{xunit=15pt,yunit=15pt}\psgrid[subgriddiv=1,griddots=4,
gridlabels=4pt](0,0)(25,4)

\psline(0,0)(1,0)(2,1)(2,0)(3,1)(5,3)(6,2)(7,1)(9,3)(9,2)(10,3)(10,2)(11,1)(12,0)(13,0)(16,3)(17,3)(18,2)(18,1)(19,1)
\psline(19,1)(20,2)(22,0)(24,2)(24,1)(25,0)

\pscircle*(0,0){0.06}\pscircle*(1,0){0.06}\pscircle*(2,1){0.06}\pscircle*(2,0){0.06}
\pscircle*(3,1){0.06}\pscircle*(4,2){0.06}\pscircle*(5,3){0.06}\pscircle*(6,2){0.06}

\pscircle*(7,1){0.06}\pscircle*(8,2){0.06}\pscircle*(9,2){0.06}
\pscircle*(9,3){0.06}\pscircle*(10,2){0.06}\pscircle*(10,3){0.06}\pscircle*(11,1){0.06}
\pscircle*(12,0){0.06}\pscircle*(13,0){0.06}\pscircle*(14,1){0.06}
\pscircle*(15,2){0.06}\pscircle*(16,3){0.06}\pscircle*(17,3){0.06}\pscircle*(18,2){0.06}
\pscircle*(18,1){0.06}\pscircle*(19,1){0.06}\pscircle*(20,2){0.06}
\pscircle*(21,1){0.06}\pscircle*(22,0){0.06}\pscircle*(23,1){0.06}\pscircle*(24,2){0.06}
\pscircle*(24,1){0.06}\pscircle*(25,0){0.06}

\end{pspicture}
\end{center}\vskip0.2cm

\caption{\small A G-Motzkin path of length $25$ with two $\mathbf{ud}$-peaks, four $\mathbf{uv}$-peaks, two $\mathbf{du}$-valleys, two $\mathbf{vu}$-valleys and two $\mathbf{uvu}$ strings.}

\end{figure}
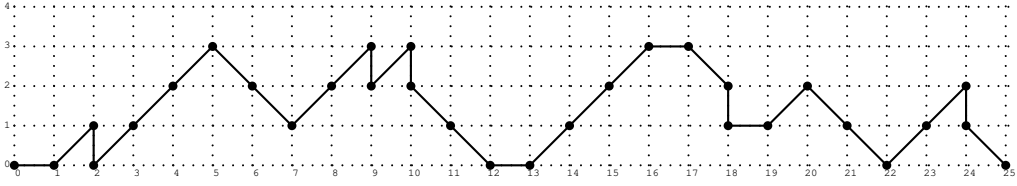

An {\it $(a,b,c)$-G-Motzkin path} is a weighted G-Motzkin path $\mathbf{P}$ such that the $\mathbf{u}$-steps, $\mathbf{h}$-steps, $\mathbf{v}$-steps and $\mathbf{d}$-steps of $\mathbf{P}$ are weighted respectively by $1, a, b$ and $c$. The {\it weight} of $\mathbf{P}$, denoted by $w(\mathbf{P})$, is the product of the weight of each step of $\mathbf{P}$. For example, $w(\mathbf{uhuduuvvdhh})=a^3b^2c^2$. The {\em weight} of a subset $\mathcal{A}$ of the set of weighted G-Motzkin paths, denoted by $w(\mathcal{A})$, is the sum of the total weights of all paths in $\mathcal{A}$. Denoted by $G_n(a, b, c)$ the weight of the set $\mathcal{G}_{n}(a, b, c)$ of all $(a,b,c)$-G-Motzkin paths of length $n$. Recently Sun et al. \cite{SunZhao} have derived the generating function of $G_n(a, b, c)$ as follows.
\begin{eqnarray*}
G(a,b,c; x)=\sum_{n=0}^{\infty}G_n(a, b, c)x^n=\frac{1-ax-\sqrt{(1-ax)^2-4x(b+cx)}}{2x(b+cx)}=\frac{1}{1-ax}C\Big(\frac{x(b+cx)}{(1-ax)^2}\Big),
\end{eqnarray*}
where
\begin{eqnarray}\label{eqn 1.1}
C(x)=\sum_{n=0}^{\infty}C_nx^n=\frac{1-\sqrt{1-4x}}{2x}
\end{eqnarray}
is the generating function for the well-known Catalan number $C_n=\frac{1}{n+1}\binom{2n}{n}$.

A {\it Dyck path} of length $2n$ is a G-Motzkin path of length $2n$ with no $\mathbf{h}$-steps or $\mathbf{v}$-steps. A {\it Motzkin path} of length $n$ is a G-Motzkin path of length $n$ with no $\mathbf{v}$-steps. An {\it $(a,b)$-Dyck path} is a weighted Dyck path with $\mathbf{u}$-steps weighted by $1$, $\mathbf{d}$-steps in $\mathbf{ud}$-peaks weighted by $a$ and other $\mathbf{d}$-steps weighted by $b$. An {\it $(a,b)$-Motzkin path} of length $n$ is an $(a,0,b)$-G-Motzkin path of length $n$.
A {\it Schr\"{o}der path} of length $2n$ is a path from $(0, 0)$ to $(2n, 0)$ in the first quadrant of the XOY-plane that consists of up steps $\mathbf{u}=(1, 1)$, down steps $\mathbf{d}=(1, -1)$ and horizontal steps $\mathbf{H}=(2, 0)$. A {\it little Schr\"{o}der path} is a Schr\"{o}der path with no $\mathbf{H}$-steps on the $x$-axis. An {\it $(a,b)$-Schr\"{o}der path} is a weighted Schr\"{o}der path such that the $\mathbf{u}$-steps, $\mathbf{H}$-steps, and $\mathbf{d}$-steps are weighted respectively by $1, a$ and $b$.
An {\it $(a,b)$-little Schr\"{o}der path} is an $(a,b)$-Schr\"{o}der path with no $\mathbf{H}$-steps on the $x$-axis.

Let $\mathcal{C}_n(a,b), \mathcal{M}_n(a,b)$ and $\mathcal{S}_n(a,b)$ be respectively the sets of $(a,b)$-Dyck paths of length $2n$, $(a,b)$-Motzkin paths of length $n$ and $(a,b)$-Schr\"{o}der paths of length $2n$. Let $C_n(a,b), M_n(a,b)$ and $S_n(a,b)$ be their weights with $C_0(a,b)=M_0(a,b)=S_0(a,b)=1$ respectively.
It is not difficult to deduce that \cite{ChenPan}
\begin{eqnarray*}
C_n(a,b)  \hskip-.22cm &=&\hskip-.22cm  \sum_{k=1}^{n}\frac{1}{n}\binom{n}{k-1}\binom{n}{k}a^{k}b^{n-k},  \\
M_n(a,b)  \hskip-.22cm &=&\hskip-.22cm  \sum_{k=0}^{n}\binom{n}{2k}C_ka^{n-2k}b^{k},                  \\
S_n(a,b)  \hskip-.22cm &=&\hskip-.22cm  \sum_{k=0}^{n}\binom{n+k}{2k}C_ka^{n-k}b^{k},
\end{eqnarray*}
and their generating functions
\begin{eqnarray*}
C(a,b; x)\hskip-.22cm &=&\hskip-.22cm \sum_{n=0}^{\infty}C_n(a,b)x^n=\frac{1-(a-b)x-\sqrt{(1-(a-b)x)^2-4bx}}{2bx}, \\
M(a,b; x)\hskip-.22cm &=&\hskip-.22cm \sum_{n=0}^{\infty}M_n(a,b)x^n=\frac{1-ax-\sqrt{(1-ax)^2-4bx^2}}{2bx^2}, \\
S(a,b; x)\hskip-.22cm &=&\hskip-.22cm \sum_{n=0}^{\infty}S_n(a,b)x^n=\frac{1-ax-\sqrt{(1-ax)^2-4bx}}{2bx}.
\end{eqnarray*}
There are closely relations between $C_n(a,b), M_n(a,b)$ and $S_n(a,b)$. Exactly, Chen and Pan \cite{ChenPan} derived the following equivalent relations
\begin{eqnarray}\label{eqn 1.2}
S_n(a,b)=C_n(a+b,b)=(a+b)M_{n-1}(a+2b,(a+b)b)
\end{eqnarray}
for $n\geq 1$ and provided some combinatorial proofs. When $a=b=1$, $C_n(a,b), M_n(a,b)$ and $S_n(a,b)$ reduce to the Catalan number $C_n$ \cite{Deutsch99, Stanley}, the Motzkin number $M_n$ \cite{Aigner, BarPinSpr, DonShap} and the large Schr\"{o}der number $S_n$ \cite{BoninShap, Sloane} respectively. Write $C(x), M(x), S(x)$ instead of $C(a,b; x), M(a,b; x)$, $S(a,b; x)$ respectively when $a=b=c=1$. In the literature, there are many papers dealing with $(a,b)$-Motzkin paths. For examples, Chen and Wang \cite{ChenWang} explored the connection between noncrossing linked partitions and $(3,2)$-Motzkin paths, established a one-to-one correspondence between the set of noncrossing linked partitions of $\{1, \dots, n+1\}$ and the set of large $(3,2)$-Motzkin paths of length $n$, which leads to a simple explanation of the well-known relation between the large and the little Schr\"{o}der numbers. Yan \cite{Yan07} found a bijective proof between the set of restricted $(3,2)$-Motzkin paths of length $n$ and the set of the Schr\"{o}der paths of length $2n$.

Recently, Sun et al. \cite{SunZhao} have studied the weighted G-Motzkin paths, counted the number of G-Motzkin paths of length $n$ with given number of $\mathbf{z}$-steps, enumerated the statistics ``number of $\mathbf{z}$-steps" at given level and discussed the statistics ``number of $\mathbf{z}_1\mathbf{z}_2$-steps" in G-Motzkin paths for $\mathbf{z}, \mathbf{z}_1, \mathbf{z}_2\in \{\mathbf{u}, \mathbf{h}, \mathbf{v}, \mathbf{d}\}$. Other related lattice paths with various steps including vertical steps permitted have been considered by \cite{Dziem-A, Dziem-B, Dziem-C, YanZhang}.

In the present paper we concentrate on the $\mathbf{uvu}$-avoiding G-Motzkin paths, that is, the G-Motzkin paths with no consecutive $\mathbf{uvu}$ strings. Precisely, in the next section we first give a bijection between the set of $\mathbf{uvu}$-avoiding $(a,b,b^2)$-G-Motzkin paths of length $n$ and the set of $(a,b)$-Schr\"{o}der paths as well as the set of $(a+b,b)$-Dyck paths of length $2n$. In the third section we mainly provide two bijections, one is between the set of $\{\mathbf{uvu, uu}\}$-avoiding $(a,b,b^2)$-G-Motzkin paths of length $n$ and the set of $(a+b,ab)$-Motzkin paths of length $n$, the other is between the set of $\{\mathbf{uvu,uu}\}$-avoiding $(a,b,b^2)$-G-Motzkin paths of length $n+1$ beginning with an $\mathbf{h}$-step weighted by $a$ and the set of $(a,b)$-Dyck paths of length $2n+2$. In the final section we focus on the enumeration of statistics ``number of $\mathbf{z}$-steps" for $\mathbf{z}\in \{\mathbf{u}, \mathbf{h}, \mathbf{v}, \mathbf{d}\}$ and ``number of points" at given level in $\mathbf{uvu}$-avoiding G-Motzkin paths. These counting results are linked with Riordan arrays. \vskip0.2cm

\section{A bijection between $\mathbf{uvu}$-avoiding $(a,b,b^2)$-G-Motzkin paths and $(a,b)$-Schr\"{o}der paths }

In this section, we first consider the $\mathbf{uvu}$-avoiding $(a,b,c)$-G-Motzkin paths which involve some classical structures as special cases, and then present a bijection between the set of $\mathbf{uvu}$-avoiding $(a,b,b^2)$-G-Motzkin paths of length $n$ and the set of $(a,b)$-Schr\"{o}der paths as well as the set of $(a+b,b)$-Dyck paths of length $2n$. When $a=b=c=1$, the $(a,b,c)$-G-Motzkin paths reduce to the common G-Motzkin paths.

Denoted by $G_n^{\mathbf{uvu}}(a, b, c)$ the weight of the set $\mathcal{G}_{n}^{\mathbf{uvu}}(a, b, c)$ of all
$\mathbf{uvu}$-avoiding $(a,b,c)$-G-Motzkin paths of length $n$. Set $\mathcal{G}^{\mathbf{uvu}}(a, b, c)=\bigcup_{n\geq 0}\mathcal{G}_n^{\mathbf{uvu}}(a, b, c)$. When $a=b=c=1$, we write
$\mathcal{G}^{\mathbf{uvu}}=\mathcal{G}^{\mathbf{uvu}}(1, 1, 1), \mathcal{G}_{n}^{\mathbf{uvu}}=\mathcal{G}_{n}^{\mathbf{uvu}}(1, 1, 1), G_n^{\mathbf{uvu}}=G_n^{\mathbf{uvu}}(1,1,1)$ for short.

Let $G^{\mathbf{uvu}}(a,b,c; x)=\sum_{n=0}^{\infty}G_n^{\mathbf{uvu}}(a, b, c)x^n$ be the generating function for the $\mathbf{uvu}$-avoiding $(a,b,c)$-G-Motzkin paths. According to the method of the first return decomposition \cite{Deutsch99}, a $\mathbf{uvu}$-avoiding $(a,b,c)$-G-Motzkin path $\mathbf{P}$ can be decomposed as one of the following six forms:
$$\mathbf{P}=\varepsilon, \ \mathbf{P}=\mathbf{uv}_b, \ \mathbf{P}=\mathbf{h}_a\mathbf{Q}_1,  \ \mathbf{P}=\mathbf{uv}_b\mathbf{h}_a\mathbf{Q}_1, \ \mathbf{P}=\mathbf{u}\mathbf{Q}_3\mathbf{v}_b\mathbf{Q}_2 \ \mbox{or}\  \mathbf{P}=\mathbf{u}\mathbf{Q}_1\mathbf{d}_c\mathbf{Q}_2, $$
where $\mathbf{x}_t$ denotes the $\mathbf{x}$-steps with weight $t$, $\mathbf{Q}_1$ and $\mathbf{Q}_2 $ are (possibly empty) $\mathbf{uvu}$-avoiding $(a,b,c)$-G-Motzkin paths, and $\mathbf{Q}_3 $ is an nonempty $\mathbf{uvu}$-avoiding $(a,b,c)$-G-Motzkin paths. Then we get the relation
\begin{eqnarray}
G^{\mathbf{uvu}}(a,b,c; x) \hskip-.22cm &=&\hskip-.22cm  1+bx+axG^{\mathbf{uvu}}(a,b,c; x)+abx^2G^{\mathbf{uvu}}(a,b,c; x)         \nonumber \\
                           \hskip-.22cm & &\hskip-.22cm  \ \ +\ bx(G^{\mathbf{uvu}}(a,b,c; x)-1)G^{\mathbf{uvu}}(a,b,c; x)+cx^2G^{\mathbf{uvu}}(a,b,c; x)^2       \nonumber\\
                           \hskip-.22cm &=&\hskip-.22cm  1+bx+(a-b+abx)xG^{\mathbf{uvu}}(a,b,c; x)+(b+cx)xG^{\mathbf{uvu}}(a,b,c; x)^2.                           \label{eqn 2.1}
\end{eqnarray}
Solve this, we have
\begin{eqnarray}\label{eqn 2.2}
G^{\mathbf{uvu}}(a,b,c; x) \hskip-.22cm &=&\hskip-.22cm \frac{(1-ax)(1+bx)-\sqrt{(1-ax)^2(1+bx)^2-4x(1+bx)(b+cx)}}{2x(b+cx)} \nonumber \\
                          \hskip-.22cm &=&\hskip-.22cm  \frac{1}{1-ax}C\Big(\frac{x(b+cx)}{(1-ax)^2(1+bx)}\Big).
\end{eqnarray}
By (\ref{eqn 1.1}), taking the coefficient of $x^n$ in $G^{\mathbf{uvu}}(a,b,c; x)$, we derive that
\begin{proposition}
For any integer $n\geq 0$, there holds
\begin{eqnarray*}
G_n^{\mathbf{uvu}}(a,b,c) \hskip-.22cm &=&\hskip-.22cm \sum_{k=0}^{n}\sum_{j=0}^{k}\sum_{\ell=0}^{n-k-j}(-1)^{\ell}\binom{k}{j}\binom{k+\ell-1}{\ell}\binom{n+k-j-\ell}{2k}C_ka^{n-k-j-\ell}b^{k+\ell-j}c^{j}  \\
                          \hskip-.22cm &=&\hskip-.22cm \sum_{k=0}^{n}\sum_{j=0}^{k}\sum_{\ell=0}^{n-k-j}(-1)^{n-k-j-\ell}\binom{k}{j}\binom{2k+\ell}{\ell}\binom{n-j-\ell-1}{n-k-j-\ell}C_ka^{\ell}b^{n-2j-\ell}c^{j}.
\end{eqnarray*}
\end{proposition}

When $(a,b,c)$ is specialized, $G^{\mathbf{uvu}}(a,b,c; x)$ and $G_n^{\mathbf{uvu}}(a,b,c)$ reduce to some well-known generating functions and classical combinatorial sequences involving the Catalan numbers $C_n$, Motzkin numbers $M_n$, Schr\"{o}der numbers $S_n$, $(a+b,b)$-Catalan number $C_n(a+b,b)$, $(a,b)$-Motzkin number $M_n(a,b)$ and $(a,b)$-Schr\"{o}der number $S_n(a,b)$. See Table 2.1 for example.
\begin{center}
\begin{eqnarray*}
\begin{array}{c|c|c|c}\hline
 (a, b, c)   & G^{\mathbf{uvu}}(a,b,c; x)               &     G_n^{\mathbf{uvu}}(a,b,c)      &       Senquences                         \\[5pt]\hline
 (0, 1, 1)   & C(x)=\frac{1-\sqrt{1-4x}}{2x}            &     C_n                            &       \mbox{\cite [A000108]{Sloane} }    \\[5pt]\hline
 (1, 0, 1)   & M(x)=\frac{1-x-\sqrt{1-2x-3x^2}}{2x^2}   &     M_n                            &       \mbox{\cite [A001006]{Sloane} }     \\[5pt]\hline
 (1, 1, 1)   & S(x)=\frac{1-x-\sqrt{1-6x+x^2}}{2x}      &     S_n                            &       \mbox{\cite [A006318]{Sloane} }      \\[5pt]\hline
 (1, 0, 2)   & \frac{1-x-\sqrt{1-2x-7x^2}}{2x}          &     a_n                            &       \mbox{\cite [A025235]{Sloane} }       \\[5pt]\hline
 (-3, 4, 16) & \frac{1+3x-\sqrt{1-10x+9x^2}}{8x}        &     a_n                            &       \mbox{\cite [A059231]{Sloane} }      \\[5pt]\hline
 (a, 0, b)   & \frac{1-ax-\sqrt{(1-ax)^2-4bx^2}}{2bx^2} &    M_n(a,b)                        &          \\[5pt]\hline
 (a, b, b^2) & \frac{1-ax-\sqrt{(1-ax)^2-4bx}}{2bx}     &    C_n(a+b,b)\ or\ S_n(a,b)        &          \\[5pt]\hline
\end{array}
\end{eqnarray*}
Table 2.1. The specializations of $G^{\mathbf{uvu}}(a,b,c; x)$ and $G_n^{\mathbf{uvu}}(a,b,c)$.
\end{center}

\begin{theorem}\label{theom 2.1.1}
For any integer $n\geq 0$, there exists a bijection $\sigma$ between $\mathcal{G}_n^{\mathbf{uvu}}(a,b,b^2)$ and $\mathcal{S}_{n}(a,b)$.
\end{theorem}
\pf Given any $\mathbf{Q}\in \mathcal{G}_n^{\mathbf{uvu}}(a,b,b^2)$ for $n\geq 0$, when $n=0, 1$, we define
$$\sigma(\varepsilon)=\varepsilon, \ \sigma(\mathbf{h}_a)=\mathbf{H}_a,\ \sigma(\mathbf{uv}_b)=\mathbf{ud}_b.$$

For $n\geq 2$, $\mathbf{Q}$ is $\mathbf{uvu}$-avoiding, there are five cases to be considered to define $\sigma(\mathbf{Q})$ recursively.

\subsection*{Case 1.} When $\mathbf{Q}=\mathbf{h}_a\mathbf{Q}'$ with $\mathbf{Q}'\in \mathcal{G}_{n-1}^{\mathbf{uvu}}(a,b,b^2)$, we define $\sigma(\mathbf{Q})=\mathbf{H}_a\sigma(\mathbf{Q}')$.

\subsection*{Case 2.} When $\mathbf{Q}=\mathbf{uv}_b\mathbf{h}_a\mathbf{Q}'$ with $\mathbf{Q}'\in \mathcal{G}_{n-2}^{\mathbf{uvu}}(a,b,b^2)$, we define $\sigma(\mathbf{Q})=\mathbf{ud}_b\mathbf{H}_a\sigma(\mathbf{Q}')$.

\subsection*{Case 3.} When $\mathbf{Q}=\mathbf{u}^i\mathbf{Q}''\mathbf{v}_b^i\mathbf{Q}'$ such that $\mathbf{Q}''\in \mathcal{G}_{k}^{\mathbf{uvu}}(a,b,b^2)$ is primitive and ending with a $\mathbf{d}$-step and $\mathbf{Q}'\in \mathcal{G}_{n-k-i}^{\mathbf{uvu}}(a,b,b^2)$ for certain $2\leq k\leq n-i$ and $1\leq i\leq n-2$, we define
\begin{eqnarray*}
\sigma(\mathbf{Q})=\left\{
\begin{array}{rl}
 (\mathbf{ud}_b\mathbf{u})^{j-1}\mathbf{ud}_b\sigma(\mathbf{Q}'')\mathbf{d}_b^{j-1}\sigma(\mathbf{Q}'),  &  \mbox{if}\ i=2j-1, \\[5pt]
 \mathbf{u}(\mathbf{ud}_b\mathbf{u})^{j-1}\mathbf{ud}_b\sigma(\mathbf{Q}'')\mathbf{d}_b^{j}\sigma(\mathbf{Q}'),  & \mbox{if}\ i=2j.
\end{array}\right.
\end{eqnarray*}
Note that $\sigma(\mathbf{Q}'')$ is always primitive in this case.

\subsection*{Case 4.} When $\mathbf{Q}=\mathbf{u}^i\mathbf{Q}''\mathbf{v}_b^i\mathbf{Q}'$ such that $\mathbf{Q}''\in \mathcal{G}_{k}^{\mathbf{uvu}}(a,b,b^2)$ is not primitive and $\mathbf{Q}'\in \mathcal{G}_{n-k-i}^{\mathbf{uvu}}(a,b,b^2)$ for certain $0\leq k\leq n-i$ and $1\leq i\leq n$, we define
\begin{eqnarray*}
\sigma(\mathbf{Q})=\left\{
\begin{array}{rl}
 \mathbf{u}(\mathbf{ud}_b\mathbf{u})^{j-1}\sigma(\mathbf{Q}'')\mathbf{d}_b^{j}\sigma(\mathbf{Q}'),  &  \mbox{if}\ i=2j-1, \\[5pt]
 (\mathbf{ud}_b\mathbf{u})^{j}\sigma(\mathbf{Q}'')\mathbf{d}_b^{j}\sigma(\mathbf{Q}'),  & \mbox{if}\ i=2j.
\end{array}\right.
\end{eqnarray*}
Note that in this case $\sigma(\mathbf{Q}'')$ is always primitive if $\mathbf{Q}''\neq \varepsilon$.

\subsection*{Case 5.} When $\mathbf{Q}=\mathbf{u}\mathbf{Q}''\mathbf{d}_{b^2}\mathbf{Q}'$ with $\mathbf{Q}''\in \mathcal{G}_{k}^{\mathbf{uvu}}(a,b,b^2)$ and $\mathbf{Q}'\in \mathcal{G}_{n-k-2}^{\mathbf{uvu}}(a,b,b^2)$ for certain $0\leq k\leq n-2$, we define $\sigma(\mathbf{Q})=\mathbf{uu}\sigma(\mathbf{Q}'')\mathbf{d}_b\mathbf{d}_b\sigma(\mathbf{Q}')$.

Conversely, the inverse procedure can be handled as follows. Given any $\mathbf{P}\in \mathcal{S}_n(a,b)$ for $n\geq 0$, when $n=0, 1$, we define
$$\sigma^{-1}(\varepsilon)=\varepsilon, \ \sigma^{-1}(\mathbf{H}_a)=\mathbf{h}_a,\ \sigma^{-1}(\mathbf{ud}_b)=\mathbf{uv}_b.$$

For $n\geq 2$, there are five cases to be considered to define $\sigma^{-1}(\mathbf{P})$ recursively.

\subsection*{Case 1.} When $\mathbf{P}=\mathbf{H}_a\mathbf{P}'$ with $\mathbf{P}'\in \mathcal{S}_{n-1}(a,b)$, we define $\sigma^{-1}(\mathbf{P})=\mathbf{h}_a\sigma^{-1}(\mathbf{P}')$.

\subsection*{Case 2.} When $\mathbf{P}=\mathbf{ud}_b\mathbf{H}_a\mathbf{P}'$ with $\mathbf{P}'\in \mathcal{S}_{n-2}(a,b)$, we define $\sigma^{-1}(\mathbf{P})=\mathbf{uv}_b\mathbf{h}_a\sigma^{-1}(\mathbf{P}')$.

\subsection*{Case 3.} When $\mathbf{P}=\mathbf{u}\mathbf{d}_b\mathbf{P}''\mathbf{P}'$ such that $\mathbf{P}''\in \mathcal{S}_{k}(a,b)$ is primitive and $\mathbf{P}'\in \mathcal{S}_{n-k-1}(a,b)$ for certain $1\leq k\leq n-1$, we define $\sigma^{-1}(\mathbf{P})=\mathbf{u}\sigma^{-1}(\mathbf{P}'')\mathbf{v}_b\sigma^{-1}(\mathbf{P}')$.

\subsection*{Case 4.} When $\mathbf{P}=\mathbf{u}\mathbf{P}''\mathbf{d}_{b}\mathbf{P}'$ such that $\mathbf{P}''\in \mathcal{S}_{k}(a,b)$ is not primitive and $\mathbf{P}'\in \mathcal{S}_{n-k-1}(a,b)$ for certain $1\leq k\leq n-1$, we define $\sigma^{-1}(\mathbf{P})=\mathbf{u}\sigma^{-1}(\mathbf{P}'')\mathbf{v}_{b}\sigma^{-1}(\mathbf{P}')$.

\subsection*{Case 5.} When $\mathbf{P}=\mathbf{uu}\mathbf{P}''\mathbf{d}_{b}\mathbf{d}_{b}\mathbf{P}'$ with $\mathbf{P}''\in \mathcal{S}_{k}(a,b)$ and $\mathbf{P}'\in \mathcal{S}_{n-k-2}(a,b)$ for certain $0\leq k\leq n-2$, we define $\sigma^{-1}(\mathbf{P})=\mathbf{u}\sigma^{-1}(\mathbf{P}'')\mathbf{d}_{b^2}\sigma^{-1}(\mathbf{P}')$.

It is not difficult to verify that $\sigma^{-1}\sigma=\sigma\sigma^{-1}=1$, both $\sigma$ and $\sigma^{-1}$ are two weight-keeping mappings and $\sigma^{-1}(\mathbf{P})$ is
$\mathbf{uvu}$-avoiding by induction on the length of $\mathbf{P}$. Hence, $\sigma$ is a desired bijection between $\mathcal{G}_n^{\mathbf{uvu}}(a,b,b^2)$ and $\mathcal{S}_{n}(a,b)$. This completes the proof of Theorem \ref{theom 2.1.1}. \qed\vskip0.2cm

In order to give a more intuitive view on the bijection $\sigma$, a pictorial description of $\sigma$ is presented for $\mathbf{Q}=\mathbf{uuu}\mathbf{v}_b\mathbf{h}_a\mathbf{u}\mathbf{d}_{b^2}\mathbf{v}_b\mathbf{v}_b\mathbf{h}_a
\mathbf{uuuuu}\mathbf{v}_b\mathbf{d}_{b^2}\mathbf{v}_b\mathbf{v}_b\mathbf{v}_b\mathbf{u}\mathbf{d}_{b^2}\in \mathcal{G}_{15}^{\mathbf{uvu}}(a,b,b^2)$, we have
$$\sigma(\mathbf{Q})=\mathbf{u}\mathbf{d}_b\mathbf{uu}\mathbf{d}_b\mathbf{H}_a\mathbf{uu}\mathbf{d}_b\mathbf{d}_b\mathbf{d}_b\mathbf{H}_a\mathbf{u}\mathbf{d}_b
\mathbf{uu}\mathbf{d}_b\mathbf{uuu}\mathbf{d}_b\mathbf{d}_b\mathbf{d}_b\mathbf{d}_b\mathbf{uu}\mathbf{d}_b\mathbf{d}_b \in \mathcal{S}_{15}(a,b). $$

See Figure 2 for detailed illustrations.

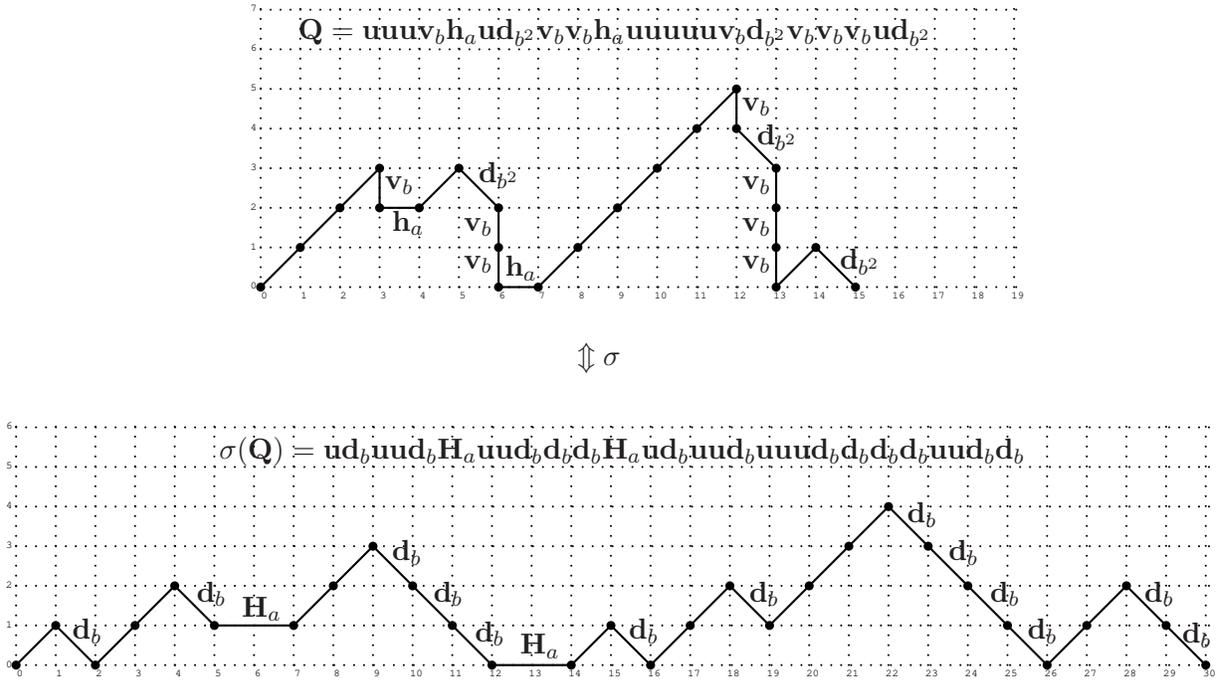
\begin{figure}[h] \setlength{\unitlength}{0.5mm}

\begin{center}
\begin{pspicture}(9,4)
\psset{xunit=15pt,yunit=15pt}\psgrid[subgriddiv=1,griddots=4,
gridlabels=4pt](0,0)(19,7)

\psline(0,0)(3,3)(3,2)(4,2)(5,3)(6,2)(6,0)(7,0)(12,5)(12,4)(13,3)(13,0)(14,1)(15,0)

\pscircle*(0,0){0.06}\pscircle*(1,1){0.06}\pscircle*(2,2){0.06}\pscircle*(3,2){0.06}
\pscircle*(3,3){0.06}\pscircle*(4,2){0.06}\pscircle*(5,3){0.06}
\pscircle*(6,2){0.06}\pscircle*(6,1){0.06}\pscircle*(6,0){0.06}\pscircle*(7,0){0.06}

\pscircle*(8,1){0.06}\pscircle*(9,2){0.06}\pscircle*(10,3){0.06}\pscircle*(11,4){0.06}
\pscircle*(12,5){0.06}\pscircle*(12,4){0.06}\pscircle*(13,3){0.06}
\pscircle*(13,2){0.06}\pscircle*(13,1){0.06}\pscircle*(13,0){0.06}\pscircle*(14,1){0.06}
 \pscircle*(15,0){0.06}

\put(1.65,1.3){$\mathbf{v}_b$}\put(1.75,0.75){$\mathbf{h}_a$}\put(2.9,1.4){$\mathbf{d}_{b^2}$}\put(2.7,0.75){$\mathbf{v}_b$}\put(2.7,0.25){$\mathbf{v}_b$}
\put(3.25,0.15){$\mathbf{h}_a$}

\put(6.6,1.9){$\mathbf{d}_{b^2}$}\put(6.4,2.35){$\mathbf{v}_b$}\put(6.4,1.3){$\mathbf{v}_b$}\put(6.4,0.75){$\mathbf{v}_b$}\put(6.4,0.25){$\mathbf{v}_b$}
\put(7.7,.25){$\mathbf{d}_{b^2}$}

\put(.5,3.3){$\mathbf{Q}=\mathbf{uuu}\mathbf{v}_b\mathbf{h}_a\mathbf{u}\mathbf{d}_{b^2}\mathbf{v}_b\mathbf{v}_b\mathbf{h}_a
\mathbf{uuuuu}\mathbf{v}_b\mathbf{d}_{b^2}\mathbf{v}_b\mathbf{v}_b\mathbf{v}_b\mathbf{u}\mathbf{d}_{b^2}$}

\end{pspicture}
\end{center}\vskip0.5cm

$\Updownarrow \sigma$
\vskip0.2cm

\begin{center}
\begin{pspicture}(18,3.3)
\psset{xunit=15pt,yunit=15pt}\psgrid[subgriddiv=1,griddots=4,
gridlabels=4pt](0,0)(30,6)

\psline(0,0)(1,1)(2,0)(4,2)(5,1)(7,1)(9,3)(12,0)(14,0)(15,1)(16,0)(18,2)(19,1)(22,4)(26,0)(28,2)(30,0)

\pscircle*(0,0){0.06}\pscircle*(1,1){0.06}\pscircle*(2,0){0.06}
\pscircle*(3,1){0.06}\pscircle*(4,2){0.06}\pscircle*(5,1){0.06}
\pscircle*(7,1){0.06}\pscircle*(8,2){0.06}
\pscircle*(9,3){0.06}\pscircle*(10,2){0.06}\pscircle*(11,1){0.06}
\pscircle*(12,0){0.06}\pscircle*(14,0){0.06}\pscircle*(15,1){0.06}
\pscircle*(16,0){0.06}\pscircle*(17,1){0.06}\pscircle*(18,2){0.06}
\pscircle*(19,1){0.06}\pscircle*(20,2){0.06}\pscircle*(21,3){0.06}\pscircle*(22,4){0.06}
\pscircle*(23,3){0.06}\pscircle*(24,2){0.06}\pscircle*(25,1){0.06}
\pscircle*(26,0){0.06}\pscircle*(27,1){0.06}\pscircle*(28,2){0.06}\pscircle*(29,1){0.06}\pscircle*(30,0){0.06}

\put(2.7,2.75){$\sigma(\mathbf{Q})=\mathbf{u}\mathbf{d}_b\mathbf{uu}\mathbf{d}_b\mathbf{H}_a
\mathbf{uu}\mathbf{d}_b\mathbf{d}_b\mathbf{d}_b\mathbf{H}_a\mathbf{u}\mathbf{d}_b
\mathbf{uu}\mathbf{d}_b\mathbf{uuu}\mathbf{d}_b\mathbf{d}_b\mathbf{d}_b\mathbf{d}_b\mathbf{uu}\mathbf{d}_b\mathbf{d}_b$}

\put(.75,.35){$\mathbf{d}_b$}\put(2.4,.85){$\mathbf{d}_b$}\put(3,.65){$\mathbf{H}_a$}\put(5,1.4){$\mathbf{d}_b$}\put(5.55,.85){$\mathbf{d}_b$}\put(6.1,.32){$\mathbf{d}_b$}
\put(6.7,.15){$\mathbf{H}_a$}\put(8.15,.35){$\mathbf{d}_b$}\put(9.75,.85){$\mathbf{d}_b$}\put(11.85,1.9){$\mathbf{d}_b$} \put(12.4,1.4){$\mathbf{d}_b$}\put(12.95,.85){$\mathbf{d}_b$}\put(13.45,.35){$\mathbf{d}_b$}\put(15,.85){$\mathbf{d}_b$}\put(15.5,.3){$\mathbf{d}_b$}

\end{pspicture}
\end{center}\vskip0.2cm

\caption{\small An example of the bijection $\sigma$ described in the proof of Theorem \ref{theom 2.1.1}. }

\end{figure}

Note that any $\mathbf{d}_{a+b}$ in $\mathbf{ud}$-peaks of $\mathbf{P}\in \mathcal{C}_n(a+b,b)$ can be regarded as
$\mathbf{d}_{a}$ and $\mathbf{d}_{b}$. If replacing all the $\mathbf{ud}_a$-peaks in $\mathbf{P}$ by $\mathbf{H}_a$ and vice versa, one can obtain a simple bijection $\phi$ between $\mathcal{C}_{n}(a+b,b)$ and $\mathcal{S}_{n}(a,b)$. For example, $$\phi(\mathbf{u}\mathbf{d}_b\mathbf{uu}\mathbf{d}_a\mathbf{u}\mathbf{d}_b\mathbf{uuu}\mathbf{d}_a\mathbf{d}_b\mathbf{d}_b\mathbf{d}_b\mathbf{u}\mathbf{d}_a)=
\mathbf{u}\mathbf{d}_b\mathbf{u}\mathbf{H}_a\mathbf{u}\mathbf{d}_b\mathbf{uu}\mathbf{H}_a\mathbf{d}_b\mathbf{d}_b\mathbf{d}_b\mathbf{H}_a.$$
Then $\phi^{-1}\sigma$ forms a bijection between $\mathcal{G}_n^{\mathbf{uvu}}(a,b,b^2)$ and $\mathcal{C}_{n}(a+b,b)$. This indicates that

\begin{corollary}\label{theom 2.1.2}
For any integer $n\geq 0$, there exists a bijection between $\mathcal{G}_n^{\mathbf{uvu}}(a,b,b^2)$ and $\mathcal{C}_{n}(a+b,b)$.
\end{corollary}

Let $\mathfrak{S}_n(a,b)$ be the set of $(a,b)$-little Schr\"{o}der paths of length $2n$ and $s_n(a,b)$ be its weight with $s_0(a,b)=1$. When $a=b=1$, $s_n=s_n(1,1)$ is the little
Schr\"{o}der number \cite[A001003]{Sloane}. Let $\hat{\mathcal{G}}_n^{\mathbf{uvu}}(a,b,b^2)$ be the set of pahts $\mathbf{Q}_1\in \mathcal{G}_n^{\mathbf{uvu}}(a,b,b^2)$ such that $\mathbf{Q}_1$ has no $\mathbf{h}$-steps on the $x$-axis. Note that under the bijection $\sigma$ any $\mathbf{h}$-step from $(i, 0)$ to $(i+1, 0)$ on the $x$-axis of $\mathbf{Q}\in \mathcal{G}_n^{\mathbf{uvu}}(a,b,b^2)$ corresponds to an $\mathbf{H}$-step from $(2i, 0)$ to $(2i+2, 0)$ on the $x$-axis of $\sigma(\mathbf{Q})\in \mathcal{S}_{n}(a,b)$ and vice versa, this signifies that the restriction of $\sigma$ to $\hat{\mathcal{G}}_n^{\mathbf{uvu}}(a,b,b^2)$ forms a bijection between $\hat{\mathcal{G}}_n^{\mathbf{uvu}}(a,b,b^2)$ and $\mathfrak{S}_n(a,b)$.

By the method of the first return decomposition \cite{Deutsch99}, the generating function $s(a,b;x)$ for $s_n(a,b)$ can be found, namely,
\begin{eqnarray*}
s(a,b; x) \hskip-.22cm &=&\hskip-.22cm \sum_{n=0}^{\infty}s_n(a,b)x^n=1+bxS(a,b; x)s(a,b; x)=\frac{1}{1-bxS(a,b; x)}  \\
          \hskip-.22cm &=&\hskip-.22cm  \frac{1+ax-\sqrt{(1-ax)^2-4bx}}{2(a+b)x}=\frac{a}{a+b}+\frac{b}{a+b}S(a,b; x),
\end{eqnarray*}
which generates
\begin{eqnarray}\label{eqn 2.3}
bS_n(a,b)=(a+b)s_n(a,b), \ \ (n\geq 1).
\end{eqnarray}
The $a=b=1$ case in (\ref{eqn 2.3}) reduces to the well-known relation $S_n=2s_n$ \cite[pp.178, 219, 256]{StanleyEC2}, whose combinatorial interpretations have been given by Shapiro and Sulanke \cite{ShapSul}, Deutsch \cite{Deutsch01}, Gu et al. \cite{GuLiMan}, Huq \cite{Huq09}, Chen and Wang \cite{ChenWang}.
Here we provide a simple combinatorial proof of (\ref{eqn 2.3}). Let $\mathbf{ud}_b\mathcal{S}_n(a,b)$ be the set of $(a,b)$-Schr\"{o}der paths of length $2n+2$ beginning with a $\mathbf{ud}_b$, and $\mathbf{ud}_b\mathfrak{S}_n(a,b)\cup \mathbf{H}_a\mathfrak{S}_n(a,b)$ be the set of $(a,b)$-little Schr\"{o}der paths of length $2n+2$ beginning with a $\mathbf{ud}_b$ or an $\mathbf{H}_a$-step. Clearly, $\mathbf{ud}_b\mathcal{S}_n(a,b)$ is counted by $bS_n(a,b)$ and $\mathbf{ud}_b\mathfrak{S}_n(a,b)\cup \mathbf{H}_a\mathfrak{S}_n(a,b)$ is counted by $(a+b)s_n(a,b)$. To prove (\ref{eqn 2.3}), it suffices to give a bijection $\vartheta$ between $\mathbf{ud}_b\mathcal{S}_n(a,b)-\mathbf{ud}_b\mathfrak{S}_n(a,b)$ and $\mathbf{H}_a\mathfrak{S}_n(a,b)$. For any $\mathbf{P}\in \mathbf{ud}_b\mathcal{S}_n(a,b)-\mathbf{ud}_b\mathfrak{S}_n(a,b)$, $\mathbf{P}$ has at least one $\mathbf{H}$-step on the $x$-axis, find the first one along the path $\mathbf{P}$, it can be partitioned into a unique form, i.e., $\mathbf{P}=\mathbf{ud}_b\mathbf{P}_1\mathbf{H}_a\mathbf{P}_2$, where $\mathbf{P}_1\in \mathfrak{S}_k(a,b)$ and $\mathbf{P}_2\in \mathcal{S}_{n-k-1}(a,b)$ for certain $0\leq k\leq n-1$. Then define
$$\vartheta(\mathbf{P})=\vartheta(\mathbf{ud}_b\mathbf{P}_1\mathbf{H}_a\mathbf{P}_2)=\mathbf{H}_a\mathbf{P}_1\mathbf{u}\mathbf{P}_2\mathbf{d}_b. $$
It is easy to verify that $\vartheta(\mathbf{P})\in \mathbf{H}_a\mathfrak{S}_n(a,b)$ with $\mathbf{u}\mathbf{P}_2\mathbf{d}_b$ being the last primitive subpath of $\vartheta(\mathbf{P})$ and that the procedure is invertible,  so $\vartheta$ is indeed a bijection between $\mathbf{ud}_b\mathcal{S}_n(a,b)-\mathbf{ud}_b\mathfrak{S}_n(a,b)$ and $\mathbf{H}_a\mathfrak{S}_n(a,b)$.  \vskip0.2cm

\section{A bijection between $\{\mathbf{uvu,uu}\}$-avoiding $(a,b,b^2)$-G-Motzkin paths and $(a+b,ab)$-Motzkin paths }

In this section we mainly provide a bijection between the set of $\{\mathbf{uvu, uu}\}$-avoiding $(a,b,b^2)$-G-Motzkin paths of length $n$ and the set of $(a+b,ab)$-Motzkin paths of length $n$. Furthermore, a bijection is also given between the set of $\{\mathbf{uvu,uu}\}$-avoiding $(a,b,b^2)$-G-Motzkin paths of length $n+1$ beginning with an $\mathbf{h}_a$-step and the set of $(a,b)$-Dyck paths of length $2n+2$.

Let $\mathcal{G}_n^{\mathbf{uvu,\tau}}(a,b,c)$ denote the set of $\{\mathbf{uvu}, \mathbf{\tau}\}$-avoiding $(a,b,c)$-G-Motzkin paths of length $n$, where $\tau$ is a single or a subset of strings on $\{\mathbf{u}, \mathbf{h}, \mathbf{v}, \mathbf{d}\}$, and denote by $G_n^{\mathbf{uvu,\tau}}(a,b,c)$ to be its weight.

\begin{theorem}\label{theom 3.1.1}
For any integer $n\geq 0$, there exists a bijection $\theta$ between $\mathcal{G}_n^{\mathbf{uvu,uu}}(a,b,b^2)$ and $\mathcal{M}_{n}(a+b,ab)$.
\end{theorem}
\pf For convenience, regard $\mathbf{h}_{a+b}$ in $\mathbf{P}\in \mathcal{M}_{n}(a+b,ab)$ as $\mathbf{h}_{a}$ and $\mathbf{h}_{b}$ because an $\mathbf{h}$-step weighted by $a+b$ is equivalent to an $\mathbf{h}$-step weighted by $a$ or $b$. Given any $\mathbf{Q}\in \mathcal{G}_n^{\mathbf{uvu, uu}}(a,b,b^2)$ for $n\geq 0$, when $n=0, 1$, we define $\theta(\varepsilon)=\mathbf{\varepsilon}, \theta(\mathbf{h}_a)=\mathbf{h}_a$ and $\theta(\mathbf{uv}_b)=\mathbf{h}_b$. When $n\geq 2$, there are five cases to be considered.

\subsection*{Case 1.} When $\mathbf{Q}=\mathbf{h}_a\mathbf{Q}_1$, where $\mathbf{Q}_1\in \mathcal{G}_k^{\mathbf{uvu, uu}}(a,b,b^2)$ for $0\leq k\leq n-1$, we define
\begin{eqnarray*}
\theta(\mathbf{Q})=\theta(\mathbf{h}_a\mathbf{Q}_1)=\mathbf{h}_a\theta(\mathbf{Q}_1).
\end{eqnarray*}

\subsection*{Case 2.} When $\mathbf{Q}=\mathbf{u}\mathbf{d}_{b^2}\mathbf{Q}_1$, where $\mathbf{Q}_1\in \mathcal{G}_k^{\mathbf{uvu, uu}}(a,b,b^2)$ for $0\leq k\leq n-2$, we define
\begin{eqnarray*}
\theta(\mathbf{Q})=\theta(\mathbf{u}\mathbf{d}_{b^2}\mathbf{Q}_1)=\mathbf{h}_{b}\mathbf{h}_{b}\theta(\mathbf{Q}_1).
\end{eqnarray*}

\subsection*{Case 3.} When $\mathbf{Q}=\mathbf{u}\mathbf{Q}_2\mathbf{d}_{b^2}\mathbf{Q}_1$, where $\mathbf{Q}_1\in \mathcal{G}_k^{\mathbf{uvu, uu}}(a,b,b^2)$ and
$\mathbf{Q}_2\in \mathcal{G}_{n-k-2}^{\mathbf{uvu, uu}}(a,b,b^2)$ is nonempty for $0\leq k\leq n-3$, $\mathbf{Q}_2$ must begin with an $\mathbf{h}$-step weighted by $a$ since $\mathbf{Q}$ is $\mathbf{uu}$-avoiding, set $\mathbf{Q}_2=\mathbf{h}_a\mathbf{Q}_2'$, we define
\begin{eqnarray*}
\theta(\mathbf{Q})=\theta(\mathbf{u}\mathbf{h}_a\mathbf{Q}_2'\mathbf{d}_{b^2}\mathbf{Q}_1)=\mathbf{h}_{b}\mathbf{u}\theta(\mathbf{Q}_2')\mathbf{d}_{ab}\theta(\mathbf{Q}_1).
\end{eqnarray*}

\subsection*{Case 4.} When $\mathbf{Q}=\mathbf{u}\mathbf{v}_{b}\mathbf{Q}_1$, where $\mathbf{Q}_1\in \mathcal{G}_k^{\mathbf{uvu, uu}}(a,b,b^2)$ for $1\leq k\leq n-1$, $\mathbf{Q}_1$ must begin with an $\mathbf{h}$-step weighted by $a$ since $\mathbf{Q}$ is $\mathbf{uvu}$-avoiding, set $\mathbf{Q}_1=\mathbf{h}_a\mathbf{Q}_1'$, we define
\begin{eqnarray*}
\theta(\mathbf{Q})=\theta(\mathbf{u}\mathbf{v}_{b}\mathbf{h}_a\mathbf{Q}_1')=\mathbf{h}_{b}\mathbf{h}_{a}\theta(\mathbf{Q}_1').
\end{eqnarray*}

\subsection*{Case 5.} When $\mathbf{Q}=\mathbf{u}\mathbf{Q}_2\mathbf{v}_{b}\mathbf{Q}_1$, where $\mathbf{Q}_1\in \mathcal{G}_k^{\mathbf{uvu, uu}}(a,b,b^2)$ and
$\mathbf{Q}_2\in \mathcal{G}_{n-k-1}^{\mathbf{uvu, uu}}(a,b,b^2)$ is nonempty for $0\leq k\leq n-2$, $\mathbf{Q}_2$ must begin with an $\mathbf{h}$-step weighted by $a$  since $\mathbf{Q}$ is $\mathbf{uu}$-avoiding, set $\mathbf{Q}_2=\mathbf{h}_a\mathbf{Q}_2'$, we define
\begin{eqnarray*}
\theta(\mathbf{Q})=\theta(\mathbf{u}\mathbf{h}_a\mathbf{Q}_2'\mathbf{v}_{b}\mathbf{Q}_1)=\mathbf{u}\theta(\mathbf{Q}_2')\mathbf{d}_{ab}\theta(\mathbf{Q}_1).
\end{eqnarray*}

Conversely, the inverse procedure can be handled as follows. Given any $\mathbf{P}\in \mathcal{M}_{n}(a+b,ab)$ for $n\geq 0$, $\theta^{-1}$ can be recursively defined as follows together with $\theta^{-1}(\varepsilon)=\varepsilon, \theta^{-1}(\mathbf{h}_a)=\mathbf{h}_a$ and $\theta^{-1}(\mathbf{h}_b)=\mathbf{uv}_b$. When $n\geq 2$,
there are five cases to be considered to define $\theta^{-1}(\mathbf{P})$ recursively.

\subsection*{Case 1.} When $\mathbf{P}=\mathbf{h}_a\mathbf{P}_1$ and $\mathbf{P}_1\in \mathcal{M}_{n-1}(a+b,ab)$, we define
\begin{eqnarray*}
\theta^{-1}(\mathbf{P})=\theta^{-1}(\mathbf{h}_a\mathbf{P}_1)=\mathbf{h}_a\theta^{-1}(\mathbf{P}_1).
\end{eqnarray*}

\subsection*{Case 2.} When $\mathbf{P}=\mathbf{h}_b\mathbf{h}_b\mathbf{P}_1$ and $\mathbf{P}_1\in \mathcal{M}_{k}(a+b,ab)$ for $0\leq k\leq n-2$, we define
\begin{eqnarray*}
\theta^{-1}(\mathbf{P})=\theta^{-1}(\mathbf{h}_b\mathbf{h}_b\mathbf{P}_1)=\mathbf{ud}_{b^2}\theta^{-1}(\mathbf{P}_1).
\end{eqnarray*}

\subsection*{Case 3.} When $\mathbf{P}=\mathbf{h}_b\mathbf{u}\mathbf{P}_2\mathbf{d}_{ab}\mathbf{P}_1$, where $\mathbf{P}_1\in \mathcal{M}_{k}(a+b,ab)$ and
$\mathbf{P}_2\in \mathcal{M}_{n-k-3}(a+b,ab)$ for $0\leq k\leq n-3$, we define
\begin{eqnarray*}
\theta^{-1}(\mathbf{P})=\theta^{-1}(\mathbf{h}_b\mathbf{u}\mathbf{P}_2\mathbf{d}_{ab}\mathbf{P}_1)=\mathbf{uh}_{a}\theta^{-1}(\mathbf{P}_2)\mathbf{d}_{b^2}\theta^{-1}(\mathbf{P}_1).
\end{eqnarray*}

\subsection*{Case 4.} When $\mathbf{P}=\mathbf{h}_b\mathbf{h}_a\mathbf{P}_1$ and $\mathbf{P}_1\in \mathcal{M}_{k}(a+b,ab)$ for $0\leq k\leq n-2$, we define
\begin{eqnarray*}
\theta^{-1}(\mathbf{P})=\theta^{-1}(\mathbf{h}_b\mathbf{h}_a\mathbf{P}_1)=\mathbf{uv}_{b}\mathbf{h}_a\theta^{-1}(\mathbf{P}_1).
\end{eqnarray*}

\subsection*{Case 5.} When $\mathbf{P}=\mathbf{u}\mathbf{P}_2\mathbf{d}_{ab}\mathbf{P}_1$, where $\mathbf{P}_1\in \mathcal{M}_{k}(a+b,ab)$, and
$\mathbf{P}_2\in \mathcal{M}_{n-k-2}(a+b,ab)$ is nonempty for $0\leq k\leq n-3$, we define
\begin{eqnarray*}
\theta^{-1}(\mathbf{P})=\theta^{-1}(\mathbf{u}\mathbf{P}_2\mathbf{d}_{ab}\mathbf{P}_1)=\mathbf{u}\mathbf{h}_{a}\theta^{-1}(\mathbf{P}_2)\mathbf{v}_{b}\theta^{-1}(\mathbf{P}_1).
\end{eqnarray*}

It can be verified that $\theta^{-1}\theta=\theta\theta^{-1}=1$, both $\theta$ and $\theta^{-1}$ are two weight-keeping mappings and $\theta^{-1}(\mathbf{P})$ is
$\{\mathbf{uvu,uu}\}$-avoiding by induction on the length of $\mathbf{P}$. Hence, $\theta$ is a desired bijection between $\mathcal{G}_n^{\mathbf{uvu,uu}}(a,b,b^2)$ and $\mathcal{M}_{n}(a+b,ab)$. This completes the proof of Theorem \ref{theom 3.1.1}. \qed\vskip0.2cm

In order to give a more intuitive view on the bijection $\theta$, a pictorial description of $\theta$ is presented for $\mathbf{Q}=\mathbf{u}\mathbf{h}_a^3\mathbf{u}\mathbf{h}_a
\mathbf{u}\mathbf{h}_a^2\mathbf{u}\mathbf{h}_a\mathbf{u}\mathbf{v}_b\mathbf{h}_a
\mathbf{v}_b^2\mathbf{u}\mathbf{d}_{b^2}^2\mathbf{v}_{b}\mathbf{h}_a^3\mathbf{u}\mathbf{h}_a^2\mathbf{u}\mathbf{h}_a\mathbf{v}_b\mathbf{u}\mathbf{v}_b\mathbf{d}_{b^2}
\mathbf{h}_a\mathbf{u}\mathbf{h}_a\mathbf{d}_{b^2} \in \mathcal{G}_{30}^{\mathbf{uvu,uu}}(a,b,b^2)$, we have
$\theta(\mathbf{Q})=\mathbf{u}\mathbf{h}_a^2\mathbf{h}_b\mathbf{uu}
\mathbf{h}_a\mathbf{u}\mathbf{h}_b\mathbf{h}_a\mathbf{d}_{ab}^2\mathbf{h}_b^2
\mathbf{d}_{ab}^2\mathbf{h}_a^3\mathbf{h}_b\mathbf{u}\mathbf{h}_a \mathbf{u}\mathbf{d}_{ab}\mathbf{h}_b\mathbf{d}_{ab}\mathbf{h}_a\mathbf{h}_b\mathbf{u}\mathbf{d}_{ab} \in \mathcal{M}_{30}(a+b,ab)$. See Figure 3 for detailed illustrations.

\begin{figure}[h] \setlength{\unitlength}{0.5mm}

\begin{center}
\begin{pspicture}(18,4.)
\psset{xunit=15pt,yunit=15pt}\psgrid[subgriddiv=1,griddots=4,
gridlabels=4pt](0,0)(30,7)

\psline(0,0)(1,1)(4,1)(5,2)(6,2)(7,3)(9,3)(10,4)(11,4)(12,5)(12,4)(13,4)(13,2)(14,3)
(15,2)(16,1)(16,0)(19,0)(20,1)(22,1)(23,2)(24,2)(24,1)(25,2)(25,1)(26,0)(27,0)(28,1)
(29,1)(30,0)

\pscircle*(0,0){0.06}\pscircle*(1,1){0.06}\pscircle*(2,1){0.06}\pscircle*(3,1){0.06}
\pscircle*(4,1){0.06}\pscircle*(5,2){0.06}\pscircle*(6,2){0.06}
\pscircle*(7,3){0.06}\pscircle*(8,3){0.06}\pscircle*(9,3){0.06}\pscircle*(10,4){0.06}

\pscircle*(11,4){0.06}\pscircle*(12,5){0.06}\pscircle*(12,4){0.06}
\pscircle*(13,4){0.06}\pscircle*(13,3){0.06}
\pscircle*(13,2){0.06}\pscircle*(14,3){0.06}\pscircle*(15,2){0.06}\pscircle*(16,1){0.06}
\pscircle*(16,0){0.06}\pscircle*(17,0){0.06}\pscircle*(18,0){0.06}\pscircle*(19,0){0.06}

\pscircle*(20,1){0.06}\pscircle*(21,1){0.06}\pscircle*(22,1){0.06}\pscircle*(23,2){0.06}
\pscircle*(24,2){0.06}\pscircle*(24,1){0.06}\pscircle*(25,2){0.06}\pscircle*(25,1){0.06}
\pscircle*(26,0){0.06}\pscircle*(27,0){0.06}\pscircle*(28,1){0.06}\pscircle*(29,1){0.06}
\pscircle*(30,0){0.06}

\put(2.7,3.3){$\mathbf{Q}=\mathbf{u}\mathbf{h}_a^3\mathbf{u}\mathbf{h}_a
\mathbf{u}\mathbf{h}_a^2\mathbf{u}\mathbf{h}_a\mathbf{u}\mathbf{v}_b\mathbf{h}_a
\mathbf{v}_b^2\mathbf{u}\mathbf{d}_{b^2}^2\mathbf{v}_{b}
\mathbf{h}_a^3\mathbf{u}\mathbf{h}_a^2\mathbf{u}\mathbf{h}_a\mathbf{v}_b\mathbf{u}\mathbf{v}_b\mathbf{d}_{b^2}
\mathbf{h}_a\mathbf{u}\mathbf{h}_a\mathbf{d}_{b^2}$}

\put(.6,.65){$\mathbf{h}_a$}\put(1.1,.65){$\mathbf{h}_a$}\put(1.65,.65){$\mathbf{h}_a$}
\put(2.75,1.15){$\mathbf{h}_a$}\put(3.8,1.7){$\mathbf{h}_a$}\put(4.3,1.7){$\mathbf{h}_a$}
\put(5.35,2.25){$\mathbf{h}_a$}\put(6.4,2.4){$\mathbf{v}_b$}

\put(6.4,1.9){$\mathbf{h}_a$}\put(6.45,1.25){$\mathbf{v}_b$}\put(6.9,1.75){$\mathbf{v}_b$}

\put(8.2,.8){$\mathbf{d}_{b^2}$}\put(7.75,1.3){$\mathbf{d}_{b^2}$}

\put(8,.15){$\mathbf{v}_b$}

\put(8.55,.15){$\mathbf{h}_a$}\put(9.05,.15){$\mathbf{h}_a$}\put(9.6,.15){$\mathbf{h}_a$}

\put(10.65,.65){$\mathbf{h}_a$}\put(11.15,.65){$\mathbf{h}_a$}\put(12.2,1.15){$\mathbf{h}_a$}

\put(12.25,.65){$\mathbf{v}_b$}\put(13.25,.65){$\mathbf{v}_b$}

\put(12.95,.1){$\mathbf{d}_{b^2}$}\put(13.75,.15){$\mathbf{h}_a$}\put(14.85,.65){$\mathbf{h}_a$}
\put(15.1,.1){$\mathbf{d}_{b^2}$}

\end{pspicture}
\end{center}\vskip0.5cm

$\Updownarrow \theta$

\begin{center}
\begin{pspicture}(18,3.9)
\psset{xunit=15pt,yunit=15pt}\psgrid[subgriddiv=1,griddots=4,
gridlabels=4pt](0,0)(30,7)

\psline(0,0)(1,1)(4,1)(5,2)(6,3)(7,3)(8,4)(9,4)(10,4)(11,3)(12,2)(14,2)
(15,1)(16,0)(19,0)(20,0)(21,1)(22,1)(23,2)(24,1)(25,1)(26,0)(27,0)(28,0)
(29,1)(30,0)

\pscircle*(0,0){0.06}\pscircle*(1,1){0.06}\pscircle*(2,1){0.06}\pscircle*(3,1){0.06}
\pscircle*(4,1){0.06}\pscircle*(5,2){0.06}\pscircle*(6,3){0.06}
\pscircle*(7,3){0.06}\pscircle*(8,4){0.06}\pscircle*(9,4){0.06}\pscircle*(10,4){0.06}

\pscircle*(11,3){0.06}\pscircle*(12,2){0.06}
\pscircle*(13,2){0.06}\pscircle*(14,2){0.06}\pscircle*(15,1){0.06}
\pscircle*(16,0){0.06}\pscircle*(17,0){0.06}\pscircle*(18,0){0.06}\pscircle*(19,0){0.06}

\pscircle*(20,0){0.06}\pscircle*(21,1){0.06}\pscircle*(22,1){0.06}\pscircle*(23,2){0.06}
\pscircle*(24,1){0.06}\pscircle*(25,1){0.06}
\pscircle*(26,0){0.06}\pscircle*(27,0){0.06}\pscircle*(28,0){0.06}\pscircle*(29,1){0.06}
\pscircle*(30,0){0.06}

\put(2.7,3.3){$\theta(\mathbf{Q})=\mathbf{u}\mathbf{h}_a^2\mathbf{h}_b\mathbf{uu}
\mathbf{h}_a\mathbf{u}\mathbf{h}_b\mathbf{h}_a\mathbf{d}_{ab}^2\mathbf{h}_b^2
\mathbf{d}_{ab}^2\mathbf{h}_a^3\mathbf{h}_b\mathbf{u}\mathbf{h}_a \mathbf{u}\mathbf{d}_{ab}\mathbf{h}_b\mathbf{d}_{ab}\mathbf{h}_a\mathbf{h}_b\mathbf{u}\mathbf{d}_{ab}$}

\put(.6,.65){$\mathbf{h}_a$}\put(1.1,.65){$\mathbf{h}_a$}\put(1.65,.65){$\mathbf{h}_b$}
\put(3.25,1.7){$\mathbf{h}_a$}\put(4.3,2.25){$\mathbf{h}_b$}\put(4.85,2.25){$\mathbf{h}_a$}

\put(5.6,1.85){$\mathbf{d}_{ab}$}\put(5.55,1.1){$\mathbf{d}_{ab}$}

\put(6.4,1.2){$\mathbf{h}_b$}\put(6.9,1.2){$\mathbf{h}_b$}

\put(7.65,.85){$\mathbf{d}_{ab}$}\put(7.6,.15){$\mathbf{d}_{ab}$}

\put(8.5,.15){$\mathbf{h}_a$}\put(9.05,.15){$\mathbf{h}_a$}\put(9.6,.15){$\mathbf{h}_a$}
\put(9.6,.15){$\mathbf{h}_a$}\put(10.1,.15){$\mathbf{h}_b$}

\put(11.15,.65){$\mathbf{h}_a$}\put(12.3,.95){$\mathbf{d}_{ab}$}\put(12.75,.6){$\mathbf{h}_b$}

\put(12.95,.1){$\mathbf{d}_{ab}$}\put(13.75,.15){$\mathbf{h}_a$}\put(14.3,.15){$\mathbf{h}_b$}
\put(15.1,.1){$\mathbf{d}_{ab}$}

\end{pspicture}
\end{center}\vskip0.2cm

\caption{\small An example of the bijection $\theta$ described in the proof of Theorem \ref{theom 3.1.1}. }

\end{figure}
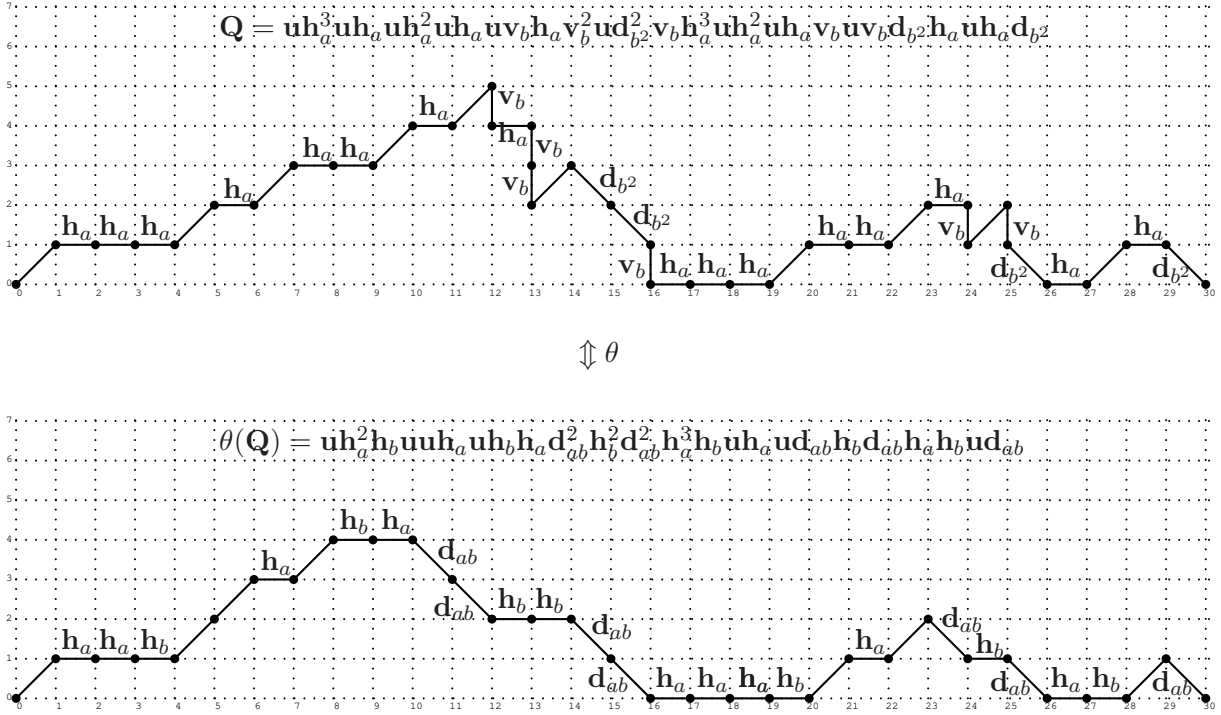

Set $\tau_1=\{\mathbf{uu},\mathbf{uh}\}$ and $\tau_2=\{\mathbf{uu},\mathbf{hu}\}$, note that any path $\mathbf{Q}\in \mathcal{G}_n^{\mathbf{uvu},\tau_1}(a,b,b^2)$ is a string of words on $\{\mathbf{h}_a, \mathbf{ud}_{b^2}, \mathbf{uv}_b\mathbf{h}_a\}$ apart from a subpath $\mathbf{uv}_b$ (if exists) at the end of $\mathbf{Q}$, and any path $\mathbf{P}\in \mathcal{G}_n^{\mathbf{uvu},\tau_2}(a,b,b^2)$ is a string of words on $\{\mathbf{u}\mathbf{h}_a^{i}\mathbf{v}_b, \mathbf{u}\mathbf{h}_a^{j}\mathbf{d}_{b^2}\}$ apart from a subpath $\mathbf{uv}_b\mathbf{h}_a^{k}$ or $\mathbf{h}_a^{\ell}$ (if exists) at the end of $\mathbf{P}$ for certain $k,\ell, j\geq 0$ and $i\geq 1$.
Let $\mathcal{B}_n$ be the set of paths $\mathbf{h}^{n}$ with $\mathbf{h}$-steps weighted by $a$ or $b$. Clearly, $\mathcal{B}_n$ is counted by $(a+b)^n$. One can verify that the restriction of $\theta$ to $\mathcal{G}_n^{\mathbf{uvu},\tau_1}(a,b,b^2)$ induces a bijection between $\mathcal{G}_n^{\mathbf{uvu},\tau_1}(a,b,b^2)$ and
$\mathcal{B}_n$. Moreover, a simple bijection $\rho$ between $\mathcal{G}_n^{\mathbf{uvu},\tau_2}(a,b,b^2)$ and $\mathcal{B}_n$ for $n\geq 0$ is defined as follows
\begin{eqnarray*}
\rho(\mathbf{P})=\left\{
\begin{array}{rl}
 \mathbf{h}_a^{n},                                           &  \mbox{if}\ \mathbf{P}=\mathbf{h}_a^{n}, \\[5pt]
 \mathbf{h}_b\mathbf{h}_a^{n-1},                             &  \mbox{if}\ \mathbf{P}=\mathbf{u}\mathbf{v}_b\mathbf{h}_a^{n-1}, \\[5pt]
 \mathbf{h}_a^{i}\mathbf{h}_b\rho(\mathbf{P}_1),             &  \mbox{if}\ \mathbf{P}=\mathbf{u}\mathbf{h}_a^{i}\mathbf{v}_b\mathbf{P}_1 \ \mbox{for}\ i\geq 1, \\[5pt]
 \mathbf{h}_b\mathbf{h}_a^{j}\mathbf{h}_b\rho(\mathbf{P}_1), &  \mbox{if}\ \mathbf{P}=\mathbf{u}\mathbf{h}_a^{j}\mathbf{d}_{b^2}P_1\ \mbox{for}\ j\geq 0. \\[5pt]
\end{array}\right.
\end{eqnarray*}
Then we have

\begin{corollary}\label{coro 3.1.2}
For any integer $n\geq 0$, there holds
\begin{eqnarray*}
G_n^{\mathbf{uvu},\tau_1}(a,b,b^2)=G_n^{\mathbf{uvu},\tau_2}(a,b,b^2)=(a+b)^n.
\end{eqnarray*}
\end{corollary}

\vskip0.2cm

Let $\mathbf{h}_a\mathcal{G}_n^{\mathbf{uvu,uu}}(a,b,b^2)$ and $\mathbf{h}_a\mathcal{M}_n(a+b,ab)$ denote the subset of paths in $\mathcal{G}_{n+1}^{\mathbf{uvu,uu}}(a,b,b^2)$ and in $\mathcal{M}_{n+1}(a+b,ab)$ respectively both starting with an $\mathbf{h}_a$-step. By (\ref{eqn 1.2}) and Theorem \ref{theom 3.1.1},
both $\mathbf{h}_a\mathcal{G}_n^{\mathbf{uvu,uu}}(a,b,b^2)$ and $\mathbf{h}_a\mathcal{M}_n(a+b,ab)$ are counted by $aM_n(a+b,ab)=C_{n+1}(a,b)$.
In order to give a bijection between $\mathbf{h}_a\mathcal{G}_n^{\mathbf{uvu,uu}}(a,b,b^2)$ and $\mathcal{C}_{n+1}(a,b)$, by Theorem \ref{theom 3.1.1}, it suffices to present
a bijection between $\mathbf{h}_a\mathcal{M}_n(a+b,ab)$ and $\mathcal{C}_{n+1}(a,b)$.

\begin{lemma}\label{lemma 3.1.3}
For any integer $n\geq 0$, there exists a bijection $\varphi$ between $\mathbf{h}_a\mathcal{M}_n(a+b,ab)$ and $\mathcal{C}_{n+1}(a,b)$.
\end{lemma}

\pf Given any $\mathbf{h}_a\mathbf{Q}\in \mathbf{h}_a\mathcal{M}_n(a+b,ab)$ for $n\geq 0$, when $n=0, 1$, we define $\varphi(\mathbf{h}_a)=\mathbf{ud}_a$, $\varphi(\mathbf{h}_a\mathbf{h}_a)=\mathbf{ud}_a\mathbf{ud}_a$ and $\varphi(\mathbf{h}_a\mathbf{h}_b)=\mathbf{uu}\mathbf{d}_a\mathbf{d}_b$. When $n\geq 2$, there are three cases to be considered.

\subsection*{Case 1.} When $\mathbf{h}_a\mathbf{Q}=\mathbf{h}_a\mathbf{Q}_1\mathbf{h}_a$, define
\begin{eqnarray*}
\varphi(\mathbf{h}_a\mathbf{Q})=\varphi(\mathbf{h}_a\mathbf{Q}_1\mathbf{h}_{a})=\varphi(\mathbf{h}_a\mathbf{Q}_1)\mathbf{u}\mathbf{d}_a.
\end{eqnarray*}

\subsection*{Case 2.} When $\mathbf{h}_a\mathbf{Q}=\mathbf{h}_a\mathbf{Q}_1\mathbf{h}_b$, define
\begin{eqnarray*}
\varphi(\mathbf{h}_a\mathbf{Q})=\varphi(\mathbf{h}_a\mathbf{Q}_1\mathbf{h}_{b})=\mathbf{u}\varphi(\mathbf{h}_a\mathbf{Q}_1)\mathbf{d}_b.
\end{eqnarray*}

\subsection*{Case 3.} When $\mathbf{h}_a\mathbf{Q}=\mathbf{h}_a\mathbf{Q}_2\mathbf{u}\mathbf{Q}_1\mathbf{d}_{ab}$, define
\begin{eqnarray*}
\varphi(\mathbf{h}_a\mathbf{Q})=\varphi(\mathbf{h}_a\mathbf{Q}_2\mathbf{u}\mathbf{Q}_1\mathbf{d}_{ab})=
\varphi(\mathbf{h}_a\mathbf{Q}_2)\mathbf{u}\varphi(\mathbf{h}_a\mathbf{Q}_1)\mathbf{d}_b.
\end{eqnarray*}

The inverse procedure can be handled similarly, which is left to the interested readers. This completes the proof of Lemma \ref{lemma 3.1.3}.  \qed \vskip0.2cm

By Theorem \ref{theom 3.1.1} and Lemma \ref{lemma 3.1.3}, it is clear that $\varphi\theta$ forms a bijection between $\mathbf{h}_a\mathcal{G}_n^{\mathbf{uvu,uu}}(a,b,b^2)$ and $\mathcal{C}_{n+1}(a,b)$, so one gets that

\begin{theorem}\label{theom 3.1.4}
For any integer $n\geq 0$, there exists a bijection between $\mathbf{h}_a\mathcal{G}_n^{\mathbf{uvu,uu}}(a,b,b^2)$ and $\mathcal{C}_{n+1}(a,b)$.
\end{theorem}\vskip0.2cm

\begin{remark}
Recall that $\phi^{-1}\sigma$ forms a bijection between $\mathcal{G}_n^{\mathbf{uvu}}(a,b,b^2)$ and $\mathcal{C}_{n}(a+b,b)$, by Lemma \ref{lemma 3.1.3}, one has that $\psi=\varphi^{-1}\phi^{-1}\sigma$ forms a bijection between $\mathcal{G}_{n}^{\mathbf{uvu}}(a,b,b^2)$ and $\mathbf{h}_{a+b}\mathcal{M}_{n-1}(a+2b,(a+b)b)$. In order to understand these bijections better, see Figure 4 for a digraph illustration. It signifies that we provide another combinatorial interpretations for the relation (\ref{eqn 1.2}).

\begin{center}
\begin{eqnarray*}
\begin{array}{rcccl}
       &          &  \mathcal{G}_n^{\mathbf{uvu}}(a,b,b^2)               &    \autoleftrightharpoons{$\sigma^{-1}$}{$\sigma$}       &   \mathcal{S}_{n}(a,b)      \\[10pt]
       &          &  \psi\downharpoonleft \upharpoonright \psi^{-1}      &                                         &    \phi^{-1} \downharpoonleft \upharpoonright\phi \\[5pt]
\mathbf{h}_{z}\mathcal{G}_{n-1}^{\mathbf{uvu,uu}}(z,b,b^2)|_{z=a+b}   &    \autoleftrightharpoons{$\theta^{-1}$}{$\theta$}    &  \mathbf{h}_{z}\mathcal{M}_{n-1}(z+b,zb)|_{z=a+b}      &    \autoleftrightharpoons{$\varphi^{-1}$}{$\varphi$}     &   \mathcal{C}_{n}(z,b)|_{z=a+b}  \\
\end{array}
\end{eqnarray*}\vskip0.2cm
Figure 4. A digraph illustration for the bijections $\sigma, \phi, \varphi, \psi$ and $\theta$.
\end{center}

\end{remark}

In order to give a more intuitive view on the bijection $\varphi\theta$, a pictorial description of $\varphi\theta$ is presented for $\mathbf{Q}=\mathbf{h}_a\mathbf{u}\mathbf{h}_a\mathbf{v}_b\mathbf{u}
\mathbf{h}_a\mathbf{d}_{b^2}\mathbf{h}_a\mathbf{u}\mathbf{h}_a\mathbf{u}\mathbf{h}_a\mathbf{u}\mathbf{d}_{b^2}\mathbf{u}\mathbf{v}_{b}\mathbf{d}_{b^2}\mathbf{v}_b \in \mathbf{h}_a\mathcal{G}_n^{\mathbf{uvu}}(a,b,b^2)$, we have
\begin{eqnarray*}
\theta(\mathbf{Q})  \hskip-.22cm &=&\hskip-.22cm        \mathbf{h}_a\mathbf{u}\mathbf{d}_{ab}\mathbf{h}_b\mathbf{u}\mathbf{d}_{ab}\mathbf{h}_a\mathbf{u}
                                                        \mathbf{h}_b\mathbf{u}\mathbf{h}_b\mathbf{h}_b\mathbf{h}_b\mathbf{d}_{ab}\mathbf{d}_{ab} \\
\varphi\theta(\mathbf{Q}) \hskip-.22cm &=&\hskip-.22cm  \mathbf{u}\mathbf{u}\mathbf{d}_{a}\mathbf{u}^2\mathbf{d}_{a}\mathbf{d}_{b}\mathbf{d}_{b}\mathbf{u}^2\mathbf{d}_{a}\mathbf{d}_{b}
                                                        \mathbf{u}\mathbf{d}_{a}\mathbf{u}\mathbf{u}^2\mathbf{d}_{a}\mathbf{d}_{b}\mathbf{u}^5\mathbf{d}_{a}\mathbf{d}_{b}^4\mathbf{d}_{b}.
\end{eqnarray*}
See Figure 5 for detailed illustrations.

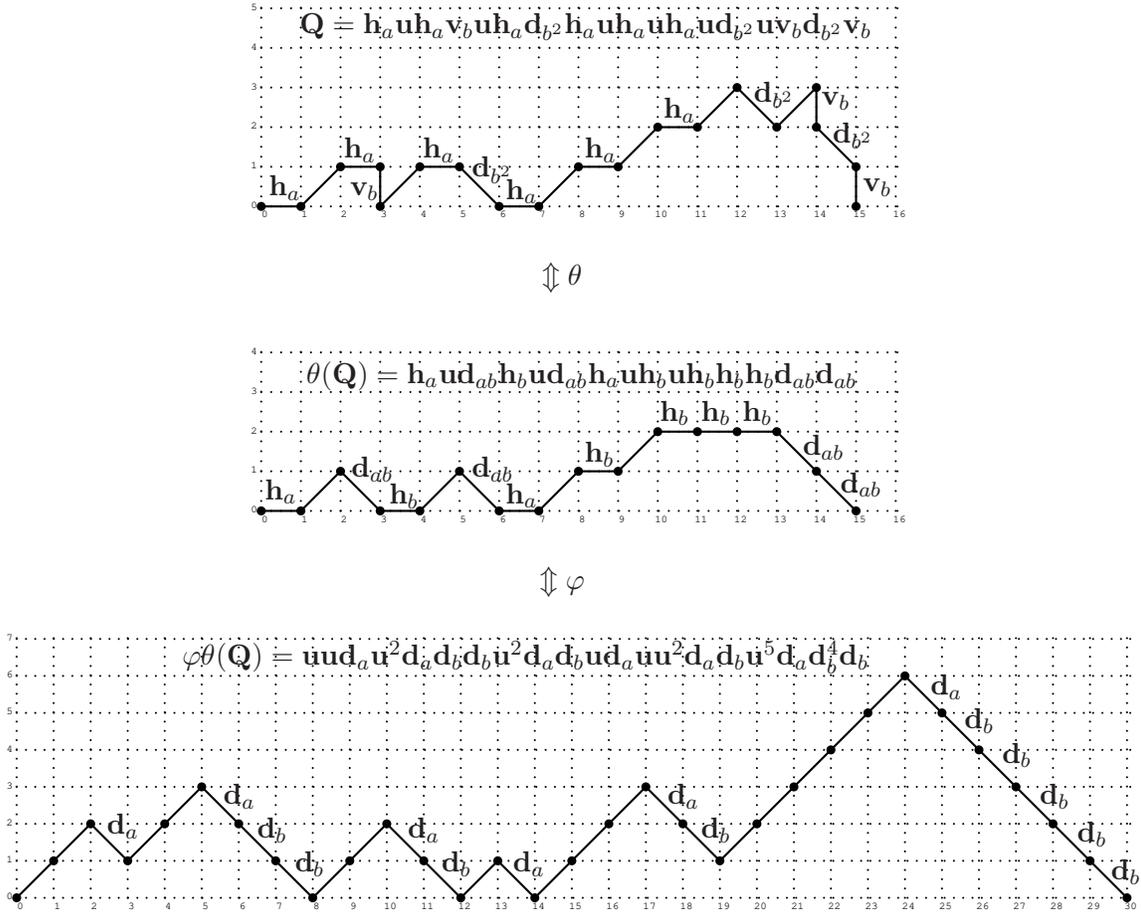
\begin{figure}[h] \setlength{\unitlength}{0.5mm}

\begin{center}
\begin{pspicture}(8,3.2)
\psset{xunit=15pt,yunit=15pt}\psgrid[subgriddiv=1,griddots=4,
gridlabels=4pt](0,0)(16,5)

\psline(0,0)(1,0)(2,1)(3,1)(3,0)(4,1)(5,1)(6,0)(7,0)(8,1)(9,1)(10,2)(11,2)(12,3)(13,2)(14,3)(14,2)(15,1)(15,0)

\pscircle*(0,0){0.06}\pscircle*(1,0){0.06}\pscircle*(2,1){0.06}\pscircle*(3,0){0.06}
\pscircle*(3,1){0.06}\pscircle*(4,1){0.06}\pscircle*(5,1){0.06}\pscircle*(6,0){0.06}

\pscircle*(7,0){0.06}\pscircle*(8,1){0.06} \pscircle*(9,1){0.06}\pscircle*(10,2){0.06}\pscircle*(11,2){0.06}
\pscircle*(12,3){0.06}\pscircle*(13,2){0.06}\pscircle*(14,2){0.06}\pscircle*(14,3){0.06}\pscircle*(15,0){0.06}
\pscircle*(15,1){0.06}

\put(.1,.15){$\mathbf{h}_a$}\put(1.1,.65){$\mathbf{h}_a$}\put(1.18,.15){$\mathbf{v}_b$}\put(2.15,.65){$\mathbf{h}_a$}
\put(2.8,.4){$\mathbf{d}_{b^2}$}\put(3.25,.1){$\mathbf{h}_a$}\put(4.3,.65){$\mathbf{h}_a$}\put(5.35,1.18){$\mathbf{h}_a$}
\put(6.55,1.4){$\mathbf{d}_{b^2}$}\put(7.45,1.35){$\mathbf{v}_b$}\put(7.6,0.85){$\mathbf{d}_{b^2}$}\put(8,.2){$\mathbf{v}_b$}

\put(.51,2.3){$\mathbf{Q}=\mathbf{h}_a\mathbf{u}\mathbf{h}_a\mathbf{v}_b\mathbf{u}
\mathbf{h}_a\mathbf{d}_{b^2}\mathbf{h}_a\mathbf{u}\mathbf{h}_a\mathbf{u}\mathbf{h}_a\mathbf{u}\mathbf{d}_{b^2}\mathbf{u}\mathbf{v}_{b}\mathbf{d}_{b^2}\mathbf{v}_b$}

\end{pspicture}
\end{center}\vskip0.5cm

$\Updownarrow \theta$

\begin{center}
\begin{pspicture}(8,2.8)
\psset{xunit=15pt,yunit=15pt}\psgrid[subgriddiv=1,griddots=4,
gridlabels=4pt](0,0)(16,4)

\psline(0,0)(1,0)(2,1)(3,0)(4,0)(5,1)(6,0)(7,0)(8,1)(9,1)(10,2)(11,2)(12,2)(13,2)(14,1)(15,0)

\pscircle*(0,0){0.06}\pscircle*(1,0){0.06}\pscircle*(2,1){0.06}\pscircle*(3,0){0.06}
\pscircle*(4,0){0.06}\pscircle*(5,1){0.06}\pscircle*(6,0){0.06}
\pscircle*(7,0){0.06}\pscircle*(8,1){0.06}\pscircle*(9,1){0.06}\pscircle*(10,2){0.06}

\pscircle*(11,2){0.06}\pscircle*(12,2){0.06}\pscircle*(13,2){0.06}\pscircle*(14,1){0.06}\pscircle*(15,0){0.06}

\put(.6,1.7){$\theta(\mathbf{Q})=\mathbf{h}_a\mathbf{u}\mathbf{d}_{ab}\mathbf{h}_b\mathbf{u}\mathbf{d}_{ab}\mathbf{h}_a\mathbf{u}
\mathbf{h}_b\mathbf{u}\mathbf{h}_b\mathbf{h}_b\mathbf{h}_b\mathbf{d}_{ab}\mathbf{d}_{ab} $}

\put(.05,.15){$\mathbf{h}_a$}\put(1.2,.45){$\mathbf{d}_{ab}$}\put(1.7,.1){$\mathbf{h}_b$}
\put(2.8,.45){$\mathbf{d}_{ab}$}\put(3.25,.1){$\mathbf{h}_a$}\put(4.3,.65){$\mathbf{h}_b$}
\put(5.3,1.2){$\mathbf{h}_b$}\put(5.85,1.2){$\mathbf{h}_b$}\put(6.4,1.2){$\mathbf{h}_b$}
\put(7.2,.75){$\mathbf{d}_{ab}$}\put(7.7,.25){$\mathbf{d}_{ab}$}

\end{pspicture}
\end{center}\vskip0.5cm

$\Updownarrow  \varphi$

\begin{center}
\begin{pspicture}(14.5,3.9)
\psset{xunit=14pt,yunit=14pt}\psgrid[subgriddiv=1,griddots=4,
gridlabels=4pt](0,0)(30,7)

\psline(0,0)(2,2)(3,1)(5,3)(8,0)(10,2)(12,0)(13,1)(14,0)(17,3)(19,1)(24,6)(30,0)

\pscircle*(0,0){0.06}\pscircle*(1,1){0.06}\pscircle*(2,2){0.06}
\pscircle*(3,1){0.06}\pscircle*(4,2){0.06}\pscircle*(5,3){0.06}
\pscircle*(6,2){0.06}\pscircle*(7,1){0.06}\pscircle*(8,0){0.06}
\pscircle*(9,1){0.06}\pscircle*(10,2){0.06}\pscircle*(11,1){0.06}
\pscircle*(12,0){0.06}\pscircle*(13,1){0.06}\pscircle*(14,0){0.06}
\pscircle*(15,1){0.06}\pscircle*(16,2){0.06}\pscircle*(17,3){0.06}
\pscircle*(18,2){0.06}\pscircle*(19,1){0.06}

\pscircle*(30,0){0.06}\pscircle*(29,1){0.06}
\pscircle*(28,2){0.06}\pscircle*(27,3){0.06}\pscircle*(26,4){0.06}
\pscircle*(25,5){0.06}\pscircle*(24,6){0.06}\pscircle*(23,5){0.06}
\pscircle*(22,4){0.06}\pscircle*(21,3){0.06}\pscircle*(20,2){0.06}

\put(2.2,3.1){$\varphi\theta(\mathbf{Q})=\mathbf{u}\mathbf{u}\mathbf{d}_{a}\mathbf{u}^2\mathbf{d}_{a}\mathbf{d}_{b}\mathbf{d}_{b}\mathbf{u}^2\mathbf{d}_{a}\mathbf{d}_{b}
\mathbf{u}\mathbf{d}_{a}\mathbf{u}\mathbf{u}^2\mathbf{d}_{a}\mathbf{d}_{b}\mathbf{u}^5\mathbf{d}_{a}\mathbf{d}_{b}^4\mathbf{d}_{b} $}

\put(1.2,.85){$\mathbf{d}_a$}\put(2.75,1.25){$\mathbf{d}_a$}\put(3.2,0.8){$\mathbf{d}_b$}\put(3.7,0.35){$\mathbf{d}_b$}
\put(5.2,0.8){$\mathbf{d}_a$}\put(5.65,0.35){$\mathbf{d}_b$}\put(6.6,0.35){$\mathbf{d}_a$}

\put(8.65,1.25){$\mathbf{d}_a$}\put(9.1,0.8){$\mathbf{d}_b$}

\put(12.15,2.7){$\mathbf{d}_a$}\put(12.6,2.25){$\mathbf{d}_b$}\put(13.1,1.8){$\mathbf{d}_b$}
\put(13.6,1.25){$\mathbf{d}_b$}\put(14.1,0.75){$\mathbf{d}_b$}\put(14.55,0.25){$\mathbf{d}_b$}

\end{pspicture}
\end{center}

\caption{\small An example of the bijection $\varphi\theta$ described in the proof of Theorem \ref{theom 3.1.4}. }

\end{figure}

\section{The statistics ``number of $\mathbf{z}$-steps" in $\mathbf{uvu}$-avoiding G-Motzkin paths }

In this section, we focus on the enumeration of statistics ``number of $\mathbf{z}$-steps" for $\mathbf{z}\in \{\mathbf{u}, \mathbf{h}, \mathbf{v}, \mathbf{d}\}$ and ``number of points" at given level in $\mathbf{uvu}$-avoiding G-Motzkin paths. Some counting results are linked with Riordan arrays.

Recall that a {\it Riordan array} \cite{ShapB, ShapGet, Sprug} is an infinite lower triangular matrix $\mathscr{D}=(d_{n,i})_{n,i \in \mathbb{N}}$ such that its $i$-th column has generating function $d(x)h(x)^i$, where $d(x)$ and $h(x)$ are formal power series with $d(0)=1$ and $h(0)=0$. That is, the general term of $\mathscr{D}$ is $d_{n,i}=[x^n]d(x)h(x)^i$, where $[x^n]$ is the coefficient operator. The matrix $\mathscr{D}$ corresponding to the pair $d(x)$ and $h(x)$ is denoted by $(d(x),h(x))$.

\subsection{ The statistics ``number of $\mathbf{u}$-steps" at level $i+1$ }
Let $U_{n, i}$ denote the total number of $\mathbf{u}$-steps at level $i+1$ in all $\mathbf{G}\in \mathcal{G}_{n+1}^{\mathbf{uvu}}$, the set of all $\mathbf{uvu}$-avoiding $(1,1,1)$-G-Motzkin paths of length $n+1$. The first values of $U_{n, i}$ are illustrated in Table 4.1.

\begin{center}
\begin{eqnarray*}
\begin{array}{c|ccccccc}\hline
n/i & 0      & 1      & 2       & 3       & 4    & 5    & 6    \\\hline
  0 & 1      &        &         &         &      &      &     \\
  1 & 5      & 1      &         &         &      &      &     \\
  2 & 25     & 9      & 1       &         &      &      &      \\
  3 & 121    & 61     & 13      & 1       &      &      &      \\
  4 & 593    & 369    & 113     & 17      &  1   &      &      \\
  5 & 2941   & 2121   & 825     & 181     &  21  &  1   &      \\
  6 & 14777  & 11881  & 5489    & 1553    &  265 &  25  & 1     \\\hline
\end{array}
\end{eqnarray*}
Table 4.1. The first values of $U_{n, i}$.
\end{center}

\begin{theorem}\label{theom 4.1.1}
For any integers $n\geq i\geq 0$, there holds
\begin{eqnarray*}
U_{n, i}=\sum_{j=0}^{n-i}(-1)^{j}\sum_{m=0}^{n-i-j}\binom{n+m+i-j+2}{n-i-j-m}C_{m}^{(2i+3)},
\end{eqnarray*}
where $C_{m}^{(k)}=\frac{k}{m+k}\binom{2m+k}{m}$. Moreover, $U_{n,i}$ is the $(n,i)$-entry of the Riordan array
$$\Big(\frac{S(x)^3}{1+x},\ xS(x)^2\Big). $$
\end{theorem}

\pf For any $\mathbf{uvu}$-avoiding G-Motzkin path $\mathbf{G}\in \mathcal{G}^{\mathbf{uvu}}$ with at least one $\mathbf{u}$-step at level $i+1$ for $i\geq 0$, given such a $\mathbf{u}$-step, marked as $\mathbf{u}^{*}$, $\mathbf{G}$ can be partitioned into one of the following three cases,
\begin{eqnarray*}
\mathbf{G}\hskip-.22cm &=&\hskip-.22cm \mathbf{G}_0\mathbf{u}\mathbf{G}_1\dots \mathbf{u}\mathbf{G}_{i}\mathbf{u}^{*}\mathbf{v}\mathbf{G}_{i+1}\mathbf{z}_1\mathbf{\bar{G}}_{1}\dots \mathbf{z}_i\mathbf{\bar{G}}_{i}, \\
\mathbf{G}\hskip-.22cm &=&\hskip-.22cm \mathbf{G}_0\mathbf{u}\mathbf{G}_1\dots \mathbf{u}\mathbf{G}_{i}\mathbf{u}^{*}\mathbf{G}_{i+1}\mathbf{v}\mathbf{\bar{G}}_{0}\mathbf{z}_1\mathbf{\bar{G}}_{1}\dots \mathbf{z}_i\mathbf{\bar{G}}_{i},\\
\mathbf{G}\hskip-.22cm &=&\hskip-.22cm \mathbf{G}_0\mathbf{u}\mathbf{G}_1\dots \mathbf{u}\mathbf{G}_{i}\mathbf{u}^{*}\mathbf{G}_{i+1}\mathbf{d}\mathbf{\bar{G}}_{0}\mathbf{z}_1\mathbf{\bar{G}}_{1}\dots \mathbf{z}_i\mathbf{\bar{G}}_{i},
\end{eqnarray*}
where $\mathbf{z}_1, \dots, \mathbf{z}_{i} \in \{\mathbf{d}, \mathbf{v}\}$, $\mathbf{G}_0, \dots, \mathbf{G}_{i}\in \mathcal{G}^{\mathbf{uvu}}$ must not end with $\mathbf{uv}$, and $\mathbf{G}_{i+1}, \mathbf{\bar{G}}_0, \dots, \mathbf{\bar{G}}_{i} \in \mathcal{G}^{\mathbf{uvu}}$ such that $\mathbf{G}_{i+1}$ in the first case is empty or beginning with an $\mathbf{h}$-step and $\mathbf{G}_{i+1}$ in the second case is nonempty. Note that each of $\mathbf{G}_0, \dots, \mathbf{G}_{i}$
is counted by the generating function $\frac{S(x)}{1+x}$, each of $\mathbf{\bar{G}}_0, \dots, \mathbf{\bar{G}}_{i}$ is counted by the generating function $S(x)$, and $\mathbf{G}_{i+1}$ is counted respectively by $1+xS(x)$, $S(x)-1$ and $S(x)$ in the first, the second and the third case. Since each of $\mathbf{u}$, $\mathbf{h}$ and $\mathbf{d}$ produces an $x$ and each $\mathbf{v}$ leads to a $1$, this makes $\mathbf{z}_1\mathbf{z}_2\dots \mathbf{z}_{i}$ generate $(1+x)^{i}$.
According to the statistics ``number of $\mathbf{u}$-steps" at level $i+1$, all $\mathbf{G}\in \mathcal{G}^{\mathbf{uvu}}$ in the three cases produce the generating functions
\begin{eqnarray*}
u_{i1}(x)\hskip-.22cm &=&\hskip-.22cm x^{i+1}\Big(\frac{S(x)}{1+x}\Big)^{i+1}\big(1+xS(x)\big)S(x)^{i}(1+x)^{i}=x^{i+1}S(x)^{2i+1}\frac{1+xS(x)}{1+x}, \\[5pt]
u_{i2}(x)\hskip-.22cm &=&\hskip-.22cm x^{i+1}\Big(\frac{S(x)}{1+x}\Big)^{i+1}(S(x)-1)S(x)^{i+1}(1+x)^{i}=x^{i+1}S(x)^{2i+2}\frac{S(x)-1}{1+x},\\[5pt]
u_{i3}(x)\hskip-.22cm &=&\hskip-.22cm x^{i+1}\Big(\frac{S(x)}{1+x}\Big)^{i+1}xS(x)^{i+2}(1+x)^{i}=x^{i+1}S(x)^{2i+3}\frac{x}{1+x},
\end{eqnarray*}
respectively. Summarizing these, by the relation $S(x)=1+xS(x)+xS(x)^2$, one has the generating function
\begin{eqnarray*}
U_i(x)    \hskip-.22cm &=&\hskip-.22cm u_{i1}(x)+u_{i2}(x)+u_{i3}(x)= x^{i+1}S(x)^{2i+3}\frac{1}{1+x},
\end{eqnarray*}
which counts the total number $U_{n,i}$ of $\mathbf{u}$-steps at level $i+1$ in all $\mathbf{G}\in \mathcal{G}_{n+1}^{\mathbf{uvu}}$.
Hence, $U_{n,i}$ is the coefficient of $x^{n+1}$ in $U_i(x)$, namely,
\begin{eqnarray*}
U_{n,i} \hskip-.22cm &=&\hskip-.22cm [x^{n+1}]U_i(x)=[x^{n+1}]x^{i+1}S(x)^{2i+3}\frac{1}{1+x}=[x^n]x^{i}S(x)^{2i}\frac{S(x)^{3}}{1+x},
\end{eqnarray*}
which implies that $U_{n,i}$ is the $(n,i)$-entry of the Riordan array
$$\Big(\frac{S(x)^3}{1+x},\ xS(x)^2\Big). $$

By the relation $S(x)=\frac{1}{1-x}C(\frac{x}{(1-x)^2})$ and the series expansion \cite{Stanley}
\begin{eqnarray*}
C(x)^{\alpha}=\sum_{m=0}^{\infty}C_m^{(\alpha)}x^m=\sum_{m=0}^{\infty}\frac{\alpha}{m+\alpha}\binom{2m+\alpha}{m}x^m,
\end{eqnarray*}
we have
\begin{eqnarray*}
U_{n,i} \hskip-.22cm &=&\hskip-.22cm [x^n]x^{i}S(x)^{2i+3}\frac{1}{1+x}=[x^{n-i}]\frac{1}{1+x}S(x)^{2i+3} \\
        \hskip-.22cm &=&\hskip-.22cm [x^{n-i}]\sum_{j=0}^{\infty}(-1)^{j}x^{j}C\Big(\frac{x}{(1-x)^2}\Big)^{2i+3}\frac{1}{(1-x)^{2i+3}}  \\
        \hskip-.22cm &=&\hskip-.22cm \sum_{j=0}^{n-i}(-1)^{j}[x^{n-i-j}]\sum_{m=0}^{\infty}C_{m}^{(2i+3)}\frac{x^m}{(1-x)^{2m+2i+3}}  \\
        \hskip-.22cm &=&\hskip-.22cm \sum_{j=0}^{n-i}(-1)^{j}[x^{n-i-j}]\sum_{m=0}^{\infty}C_{m}^{(2i+3)}
        \sum_{r=0}^{\infty}\binom{2m+2i+r+2}{r}x^{m+r} \\
        \hskip-.22cm &=&\hskip-.22cm  \sum_{j=0}^{n-i}(-1)^{j}\sum_{m=0}^{n-i-j}\binom{n+m+i-j+2}{n-i-j-m}C_{m}^{(2i+3)}.
\end{eqnarray*}
This completes the proof of Theorem \ref{theom 4.1.1}.  \qed\vskip0.2cm

\begin{remark}
One can also consider the total number $u_{n,i}$ of $\mathbf{u}$-steps at level $i+1$ in all $\mathbf{G}\in \mathcal{G}_{n+1}^{\mathbf{uvu}}$ such that $\mathbf{G}$ has no $\mathbf{h}$-steps on the $x$-axis. Note that in the proof of Theorem \ref{theom 4.1.1}, the generating function for $\mathbf{G}_0$ is $\frac{s(x)}{1+x}$ since $\mathbf{G}_0\in \mathcal{G}^{\mathbf{uvu}}$ has no $\mathbf{h}$-steps on the $x$-axis and must not end with $\mathbf{uv}$, and the generating function for $\mathbf{\bar{G}}_{i}$ is $s(x)$ since $\mathbf{\bar{G}}_{i}\in \mathcal{G}^{\mathbf{uvu}}$ has no $\mathbf{h}$-steps on the $x$-axis, where $s(x)$ is the generating function for the little Schr\"{o}der paths and $s(x)=\frac{1}{2}(1+S(x))$. Similarly, one can derive that
$u_{n,i}$ is the $(n,i)$-entry of the Riordan array
$$\Big(\frac{s(x)^2S(x)}{1+x},\ xS(x)^2\Big). $$
\end{remark} \vskip0.1cm

\subsection{ The statistics ``number of $\mathbf{v}$-steps" and ``number of $\mathbf{d}$-steps" at level $i$ }
Let $V_{n, i}$ denote the number of $\mathbf{v}$-steps at level $i$ in all $\mathbf{G}\in \mathcal{G}_{n+1}^{\mathbf{uvu}}$ and $D_{n, i}$ denote the number of $\mathbf{d}$-steps at level $i$ in all $\mathbf{G}\in \mathcal{G}_{n+2}^{\mathbf{uvu}}$. The first values of $V_{n, i}$ are illustrated in Table 4.2 and that of $D_{n, i}$ are illustrated as $U_{n,i}$ in Table 4.1.

\begin{center}
\begin{eqnarray*}
\begin{array}{c|ccccccc}\hline
n/i & 0      & 1      & 2       & 3       & 4    & 5    & 6    \\\hline
  0 & 1      &        &         &         &      &      &     \\
  1 & 4      & 1      &         &         &      &      &     \\
  2 & 20     & 8      & 1       &         &      &      &      \\
  3 & 96     & 52     & 12      & 1       &      &      &      \\
  4 & 472    & 308    & 100     & 16      & 1    &      &      \\
  5 & 2348   & 1752   & 712     & 164     & 20   & 1    &      \\
  6 & 11836  & 9760   & 4664    & 1372    & 244  & 24   & 1      \\\hline
\end{array}
\end{eqnarray*}
Table 4.2. The first values of $V_{n, i}$.
\end{center}

\begin{theorem}\label{theom 4.1.2}
For any integers $n\geq i\geq 0$, there holds
\begin{eqnarray*}
V_{n, i} \hskip-.22cm &=&\hskip-.22cm U_{n,i}-U_{n-1,i}, \\
D_{n, i} \hskip-.22cm &=&\hskip-.22cm U_{n,i}.
\end{eqnarray*}
\end{theorem}
\pf Note that any $\mathbf{u}$-step at level $i+1$ in a $\mathbf{uvu}$-avoiding G-Motzkin path $\mathbf{G}\in \mathcal{G}_{n+1}^{\mathbf{uvu}}$ has a matching step, it is a $\mathbf{v}$-step or a $\mathbf{d}$-step at level $i$, this implies that $U_{n, i}=V_{n,i}+D_{n-1,i}$. So it suffices to prove that $D_{n, i}=U_{n,i}$.

For any $\mathbf{G}\in \mathcal{G}^{\mathbf{uvu}}$ with at least one $\mathbf{d}$-step at level $i$ for $i\geq 0$, given such a $\mathbf{d}$-step, marked as $\mathbf{d}^{*}$, $\mathbf{G}$ can be partitioned uniquely into
\begin{eqnarray*}
\mathbf{G}\hskip-.22cm &=&\hskip-.22cm \mathbf{G}_0\mathbf{u}\mathbf{G}_1\dots \mathbf{u}\mathbf{G}_{i}\mathbf{u}\mathbf{G}_{i+1}\mathbf{d}^{*}\mathbf{\bar{G}}_{0}\mathbf{z}_1\mathbf{\bar{G}}_{1}\dots \mathbf{z}_i\mathbf{\bar{G}}_{i},
\end{eqnarray*}
where $\mathbf{z}_1, \dots, \mathbf{z}_{i} \in \{\mathbf{d}, \mathbf{v}\}$, $\mathbf{G}_0, \dots, \mathbf{G}_{i}\in \mathcal{G}^{\mathbf{uvu}}$ must not end with $\mathbf{uv}$, and $\mathbf{G}_{i+1}, \mathbf{\bar{G}}_0, \dots, \mathbf{\bar{G}}_{i} \in \mathcal{G}^{\mathbf{uvu}}$. Similar to the proof of Theorem \ref{theom 4.1.1},
the total number $D_{n,i}$ of $\mathbf{d}$-steps at level $i$ in all $\mathbf{G}\in \mathcal{G}_{n+2}^{\mathbf{uvu}}$ is counted by the generating function
$$u_{i3}(x)= x^{i+1}\Big(\frac{S(x)}{1+x}\Big)^{i+1} xS(x)^{i+2}(1+x)^{i}=x^{i+2}\frac{S(x)^{2i+3}}{1+x}.$$
Namely,
\begin{eqnarray*}
D_{n,i} \hskip-.22cm &=&\hskip-.22cm [x^{n+2}]u_{i3}(x)=[x^{n+2}]x^{i+2}\frac{S(x)^{2i+3}}{1+x}=[x^{n}]x^{i}S(x)^{2i}\frac{S(x)^{3}}{1+x}= U_{n,i}.
\end{eqnarray*}
This completes the proof of Theorem \ref{theom 4.1.2}.  \qed\vskip0.2cm

\begin{remark}
One can also consider the total number $v_{n,i}$ of $\mathbf{v}$-steps at level $i$ in all $\mathbf{G}\in \mathcal{G}_{n+1}^{\mathbf{uvu}}$ and $d_{n, i}$ denote the total number of $\mathbf{d}$-steps at level $i$ in all $\mathbf{G}\in \mathcal{G}_{n+2}^{\mathbf{uvu}}$ such that $\mathbf{G}$ has no $\mathbf{h}$-steps on the $x$-axis. Similarly, one can derive that $v_{n, i}= u_{n,i}-u_{n-1,i}$ and $d_{n, i}=u_{n,i}$.
\end{remark} \vskip0.1cm

\subsection{ The statistics ``number of $\mathbf{h}$-steps" at level $i$ }
Let $H_{n, i}$ denote the number of $\mathbf{h}$-steps at level $i$ in all $\mathbf{uvu}$-avoiding G-Motzkin paths $\mathbf{G}\in \mathcal{G}_{n+1}^{\mathbf{uvu}}$, the first values of $H_{n, i}$ are illustrated in Table 4.3.

\begin{center}
\begin{eqnarray*}
\begin{array}{c|ccccccc}\hline
n/i & 0      & 1      & 2       & 3       & 4    & 5    & 6    \\\hline
  0 & 1      &        &         &         &      &      &     \\
  1 & 4      & 1      &         &         &      &      &     \\
  2 & 16     & 8      & 1       &         &      &      &      \\
  3 & 68     & 48     & 12      &  1      &      &      &      \\
  4 & 304    & 264    & 96      &  16     &  1   &      &      \\
  5 & 1412   & 1408   & 652     &  160    &  20  & 1    &      \\
  6 & 6752   & 7432   & 4080    &  1296   &  240 & 24   &  1     \\\hline
\end{array}
\end{eqnarray*}
Table 4.3. The first values of $H_{n, i}$.
\end{center}

\begin{theorem}\label{theom 4.1.3}
For any integers $n\geq i\geq 0$, there holds
\begin{eqnarray*}
H_{n,i}=\sum_{m=0}^{n-i}\binom{n+m+i+1}{n-i-m}C_{m}^{(2i+2)},
\end{eqnarray*}
where $C_{m}^{(k)}=\frac{k}{m+k}\binom{2m+k}{m}$. Moreover, $H_{n,i}$ is the $(n,i)$-entry of the Riordan array
$$\big(S(x)^2,\ xS(x)^2\big). $$
\end{theorem}
\pf For any $\mathbf{G}\in \mathcal{G}^{\mathbf{uvu}}$ with at least one $\mathbf{h}$-step at level $i$ for $i\geq 0$, given such an $\mathbf{h}$-step, marked as $\mathbf{h}^{*}$, $\mathbf{G}$ can be partitioned uniquely into
\begin{eqnarray*}
\mathbf{G}\hskip-.22cm &=&\hskip-.22cm \mathbf{G}_0\mathbf{u}\mathbf{G}_1\dots \mathbf{u}\mathbf{G}_{i-1}\mathbf{u}\mathbf{G}_{i}\mathbf{h}^{*}\mathbf{\bar{G}}_{0}\mathbf{z}_1\mathbf{\bar{G}}_{1}\dots \mathbf{z}_i\mathbf{\bar{G}}_{i},
\end{eqnarray*}
where $\mathbf{z}_1, \dots, \mathbf{z}_{i} \in \{\mathbf{d}, \mathbf{v}\}$, $\mathbf{G}_0, \dots, \mathbf{G}_{i-1}\in \mathcal{G}^{\mathbf{uvu}}$ must not end with $\mathbf{uv}$, and $\mathbf{G}_{i}, \mathbf{\bar{G}}_0, \dots, \mathbf{\bar{G}}_{i} \in \mathcal{G}^{\mathbf{uvu}}$. Similar to the proof of Theorem \ref{theom 4.1.1},
the total number $H_{n,i}$ of $\mathbf{h}$-steps at level $i$ in all $\mathbf{G}\in \mathcal{G}_{n+1}^{\mathbf{uvu}}$ is counted by the generating function
$$H_{i}(x)= x^{i}\Big(\frac{S(x)}{1+x}\Big)^{i}xS(x)^{i+2}(1+x)^{i}= x^{i+1}S(x)^{2i+2}.$$
Namely,
\begin{eqnarray*}
H_{n,i} \hskip-.22cm &=&\hskip-.22cm [x^{n+1}]H_{i}(x)=[x^{n+1}]x^{i+1}S(x)^{2i+2}=[x^{n}]\big(xS(x)^{2}\big)^{i}S(x)^{2},
\end{eqnarray*}
which generates that $H_{n,i}$ is the $(n,i)$-entry of the Riordan array
$$\big(S(x)^2,\ xS(x)^2\big). $$

Similar to obtain the explicit formula for $U_{n,i}$ in Theorem \ref{theom 4.1.1}, one can derive that for $H_{n,i}$ and $H_{n,0}=S_{n+1}-S_{n}$, the detailed procedure is omitted.
This completes the proof of Theorem \ref{theom 4.1.3}.  \qed\vskip0.2cm

\begin{remark}
One can also consider the total number $h_{n,i}$ of $\mathbf{h}$-steps at level $i+1$ in all $\mathbf{G}\in \mathcal{G}_{n+2}^{\mathbf{uvu}}$ such that $\mathbf{G}$ has no $\mathbf{h}$-steps on the $x$-axis. Similarly, one can derive that $h_{n,i}$ is the $(n,i)$-entry of the Riordan array
$$\Big(s(x)^2S(x)^2,\ xS(x)^2\Big). $$
\end{remark} \vskip0.1cm

\subsection{ The statistics ``number of points" at level $i$ }
Let $P_{n, i}$ denote the number of points at level $i$ in all $\mathbf{uvu}$-avoiding G-Motzkin paths $\mathbf{G}\in \mathcal{G}_{n}^{\mathbf{uvu}}$, the first values of $P_{n, i}$ are illustrated in Table 4.4.

\begin{center}
\begin{eqnarray*}
\begin{array}{c|ccccccc}\hline
n/i & 0      & 1      & 2       & 3       & 4    & 5    & 6    \\\hline
  0 & 1      &        &         &         &      &      &     \\
  1 & 4      & 1      &         &         &      &      &     \\
  2 & 15     & 7      & 1       &         &      &      &      \\
  3 & 63     & 42     & 11      &  1      &      &      &      \\
  4 & 279    & 230    & 86      &  15     & 1    &      &      \\
  5 & 1291   & 1226   & 578     &  146    & 19   &  1   &      \\
  6 & 6159   & 6470   & 3598    &  1166   & 222  &  23   &  1     \\\hline
\end{array}
\end{eqnarray*}
Table 4.4. The first values of $P_{n, i}$.
\end{center}

\begin{theorem}\label{theom 4.1.4}
For any integers $n\geq i\geq 0$, there holds
\begin{eqnarray*}
P_{n+1,0}   \hskip-.22cm &=&\hskip-.22cm S_{n+1}+U_{n,0}+H_{n,0}, \\
P_{n+1,i+1} \hskip-.22cm &=&\hskip-.22cm \sum_{j=0}^{n-i}(-1)^{j}\sum_{m=0}^{n-i-j}\binom{n+m+i-j+3}{n-i-j-m}C_{m}^{(2i+4)} \\
            \hskip-.22cm & &\hskip-.22cm \   +\  \sum_{j=0}^{n-i-2}(-1)^{j}\sum_{m=0}^{n-i-j-2}\binom{n+m+i-j+1}{n-i-j-m-2}C_{m}^{(2i+4)}.
\end{eqnarray*}
Moreover, $P_{n+1,i+1}$ is the $(n,i)$-entry of the Riordan array
\begin{eqnarray*}
\Big(\frac{(1+x^2)S(x)^{4}}{1+x},\ xS(x)^2\Big).
\end{eqnarray*}
\end{theorem}

\pf For any $\mathbf{G}\in \mathcal{G}^{\mathbf{uvu}}$ with at least one point at level $i$ for $i\geq 0$, given such an point, marked as $(*)$, when $i=0$, there are three cases to be considered, $(1)$ $\mathbf{G}=\mathbf{G}(*)$ if the marked point lies at the end of $\mathbf{G}$, (2) $\mathbf{G}=\mathbf{\bar{G}}_1(*)\mathbf{h}\mathbf{\bar{G}}_2$ if the marked point lies at the beginning of an $\mathbf{h}$-step of $\mathbf{G}$ and (3) $\mathbf{G}=\mathbf{G}_1'(*)\mathbf{G}_2'$ if the marked point lies at the beginning of a $\mathbf{u}$-step of $\mathbf{G}$, where $\mathbf{\bar{G}}_{1}, \mathbf{\bar{G}}_{2}, \mathbf{G}_1', \mathbf{G}_2'\in \mathcal{G}^{\mathbf{uvu}}$ such that $\mathbf{G}_1'$ must not end with $\mathbf{uv}$ and $\mathbf{G}_2'$ is beginning with a $\mathbf{u}$-step. Similar to the proof of Theorem \ref{theom 4.1.1}, the three cases produce the generating functions $S(x), xS(x)^2$ and $\frac{S(x)}{1+x}xS(x)^2$ respectively. So the number $P_{n, 0}$ has the generating function
\begin{eqnarray*}
P_0(x) \hskip-.22cm &=&\hskip-.22cm  S(x)+xS(x)^2+\frac{S(x)}{1+x}xS(x)^2= S(x)+xS(x)^2+x\frac{S(x)^3}{1+x}  \\
       \hskip-.22cm &=&\hskip-.22cm  \sum_{n=0}^{\infty}S_nx^n+\sum_{n=0}^{\infty}H_{n,0}x^{n+1}+\sum_{n=0}^{\infty}U_{n,0}x^{n+1}= 1+\sum_{n=0}^{\infty}(S_{n+1}+H_{n,0}+U_{n,0})x^{n+1},
\end{eqnarray*}
which shows that $P_{n+1,0} =S_{n+1}+U_{n,0}+H_{n,0}$ with $P_{0,0}=1$.

When $i\geq 1$, there are five cases to be considered.

$\bullet$ Case (I): $\mathbf{G}=\mathbf{G}_0\mathbf{u}\mathbf{G}_1\dots \mathbf{u}\mathbf{G}_{i-1}\mathbf{u}\mathbf{G}_{i}(*)\mathbf{z}_1\mathbf{\bar{G}}_{1}\mathbf{z}_2\mathbf{\bar{G}}_{2}\dots \mathbf{z}_i\mathbf{\bar{G}}_{i}$ if the marked point lies at level $i$ and at the end of $\mathbf{G}_i$, where $\mathbf{z}_1, \dots, \mathbf{z}_{i} \in \{\mathbf{d}, \mathbf{v}\}$, $\mathbf{G}_0, \dots, \mathbf{G}_{i-1}\in \mathcal{G}^{\mathbf{uvu}}$ must not end with $\mathbf{uv}$, and $\mathbf{G}_{i}, \mathbf{\bar{G}}_1, \dots, \mathbf{\bar{G}}_{i} \in \mathcal{G}^{\mathbf{uvu}}$ such that $\mathbf{G}_{i}$ is nonempty. This case induces the generating function
$$p_{i1}(x)=x^{i}\Big(\frac{S(x)}{1+x}\Big)^{i}(S(x)-1)(1+x)^{i}S(x)^i=x^{i}S(x)^{2i}(S(x)-1).$$

$\bullet$ Case (II): $\mathbf{G}=\mathbf{G}_0\mathbf{u}\mathbf{G}_1\dots \mathbf{u}\mathbf{G}_{i-1}\mathbf{u}\mathbf{G}_{i}(*)\mathbf{h}\mathbf{G}_{i+1}\mathbf{z}_1\mathbf{\bar{G}}_{1}\mathbf{z}_2\mathbf{\bar{G}}_{2}\dots \mathbf{z}_i\mathbf{\bar{G}}_{i}$ if the marked point lies at the beginning of an $\mathbf{h}$-step at level $i$, where $\mathbf{z}_1, \dots, \mathbf{z}_{i} \in \{\mathbf{d}, \mathbf{v}\}$, $\mathbf{G}_0, \dots, \mathbf{G}_{i-1}\in \mathcal{G}^{\mathbf{uvu}}$ must not end with $\mathbf{uv}$, and $\mathbf{G}_{i}, \mathbf{G}_{i+1}, \mathbf{\bar{G}}_1, \dots, \mathbf{\bar{G}}_{i} \in \mathcal{G}^{\mathbf{uvu}}$. This case induces the generating function
$$p_{i2}(x)=x^{i}\Big(\frac{S(x)}{1+x}\Big)^{i}xS(x)^2(1+x)^{i}S(x)^i=x^{i+1}S(x)^{2i+2}.$$

$\bullet$ Case (III): $\mathbf{G}=\mathbf{G}_0\mathbf{u}\mathbf{G}_1\dots \mathbf{u}\mathbf{G}_{i-1}\mathbf{u}\mathbf{G}_{i}(*)\mathbf{G}_{i+1}\mathbf{z}_1\mathbf{\bar{G}}_{1}\mathbf{z}_2\mathbf{\bar{G}}_{2}\dots \mathbf{z}_i\mathbf{\bar{G}}_{i}$ if the marked point lies at the beginning of a $\mathbf{u}$-step at level $i+1$, where $\mathbf{z}_1, \dots, \mathbf{z}_{i} \in \{\mathbf{d}, \mathbf{v}\}$, $\mathbf{G}_0, \dots, \mathbf{G}_{i-1}, \mathbf{G}_{i}\in \mathcal{G}^{\mathbf{uvu}}$ must not end with $\mathbf{uv}$, and $\mathbf{G}_{i+1}, \mathbf{\bar{G}}_1, \dots, \mathbf{\bar{G}}_{i} \in \mathcal{G}^{\mathbf{uvu}}$ such that $\mathbf{G}_{i+1}$ is beginning with a $\mathbf{u}$-step. This case induces the generating function
$$p_{i3}(x)=x^{i}\Big(\frac{S(x)}{1+x}\Big)^{i+1}xS(x)^2(1+x)^{i}S(x)^i=x^{i+1}S(x)^{2i+3}\frac{1}{1+x}.$$

$\bullet$ Case (IV): $\mathbf{G}=\mathbf{G}_0\mathbf{u}\mathbf{G}_1\dots \mathbf{u}\mathbf{G}_{i-1}\mathbf{u}(*)\mathbf{v}\mathbf{G}_{i}\mathbf{z}_1\mathbf{\bar{G}}_{1}\mathbf{z}_2\mathbf{\bar{G}}_{2}\dots \mathbf{z}_{i-1}\mathbf{\bar{G}}_{i-1}$ if the marked point is the peak point at level $i$ of a $\mathbf{uv}$-peak, where $\mathbf{z}_1, \dots, \mathbf{z}_{i-1} \in \{\mathbf{d}, \mathbf{v}\}$, $\mathbf{G}_0, \dots, \mathbf{G}_{i-1}\in \mathcal{G}^{\mathbf{uvu}}$ must not end with $\mathbf{uv}$, and $\mathbf{G}_{i}, \mathbf{\bar{G}}_1, \dots, \mathbf{\bar{G}}_{i-1} \in \mathcal{G}^{\mathbf{uvu}}$ such that $\mathbf{G}_{i}$ must not begin with a $\mathbf{u}$-step. This case induces the generating function
$$p_{i4}(x)=x^{i}\Big(\frac{S(x)}{1+x}\Big)^{i}(1+xS(x))(1+x)^{i-1}S(x)^{i-1}=x^{i}S(x)^{2i-1}\frac{1+xS(x)}{1+x}.$$

$\bullet$ Case (V): $\mathbf{G}=\mathbf{G}_0\mathbf{u}\mathbf{G}_1\dots \mathbf{u}\mathbf{G}_{i-1}\mathbf{u}(*)\mathbf{d}\mathbf{G}_{i}\mathbf{z}_1\mathbf{\bar{G}}_{1}\mathbf{z}_2\mathbf{\bar{G}}_{2}\dots \mathbf{z}_{i-1}\mathbf{\bar{G}}_{i-1}$ if the marked point is the peak point at level $i$ of a $\mathbf{ud}$-peak, where $\mathbf{z}_1, \dots, \mathbf{z}_{i-1} \in \{\mathbf{d}, \mathbf{v}\}$, $\mathbf{G}_0, \dots, \mathbf{G}_{i-1}\in \mathcal{G}^{\mathbf{uvu}}$ must not end with $\mathbf{uv}$, and $\mathbf{G}_{i}, \mathbf{\bar{G}}_1, \dots, \mathbf{\bar{G}}_{i-1} \in \mathcal{G}^{\mathbf{uvu}}$. This case induces the generating function
$$p_{i5}(x)=x^{i}\Big(\frac{S(x)}{1+x}\Big)^{i}xS(x)(1+x)^{i-1}S(x)^{i-1}=x^{i+1}S(x)^{2i}\frac{1}{1+x}.$$

Summarizing these five cases, the total number $P_{n,i}$ of points at level $i\geq 1$ in all $\mathbf{G}\in \mathcal{G}_{n}^{\mathbf{uvu}}$ for $n\geq 1$ is counted by the  following generating function after simplification by the relation $S(x)=1+xS(x)+xS(x)^2$, i.e.,
\begin{eqnarray*}
P_{i}(x) \hskip-.22cm &=&\hskip-.22cm p_{i1}(x)+p_{i2}(x)+p_{i3}(x)+p_{i4}(x)+p_{i5}(x)= x^{i}S(x)^{2i+2}\frac{1+x^2}{1+x}.
\end{eqnarray*}

Namely, for $n\geq i\geq 0$, there holds
\begin{eqnarray*}
P_{n+1, i+1} \hskip-.22cm &=&\hskip-.22cm [x^{n+1}]P_{i+1}(x)=[x^{n+1}]x^{i+1}S(x)^{2i+4}\frac{1+x^2}{1+x} \\
             \hskip-.22cm &=&\hskip-.22cm [x^n]x^{i}S(x)^{2i}\frac{(1+x^2)S(x)^{4}}{1+x},
\end{eqnarray*}
which generates that $P_{n+1,i+1}$ is the $(n,i)$-entry of the Riordan array
$$\Big(\frac{(1+x^2)S(x)^{4}}{1+x},\ xS(x)^2\Big). $$

Similar to obtain the explicit formula for $U_{n,i}$ in Theorem \ref{theom 4.1.1}, one can derive that for $P_{n+1,i+1}$, the detailed procedure is omitted.
This completes the proof of Theorem \ref{theom 4.1.4}.  \qed\vskip0.2cm

\begin{remark}
One can also consider the total number $p_{n,i}$ of points at level $i$ in all $\mathbf{G}\in \mathcal{G}_{n}^{\mathbf{uvu}}$ such that $\mathbf{G}$ has no $\mathbf{h}$-steps on the $x$-axis. Similarly, one can derive that $p_{n,0}$ has the generating function $s(x)+\frac{xs(x)S(x)}{1+x}$ and $p_{n+1,i+1}$ is the $(n,i)$-entry of the Riordan array
$$\Big(\frac{1+x^2}{1+x}s(x)^2S(x)^2,\ xS(x)^2\Big). $$
\end{remark} \vskip0.1cm

\vskip0.5cm
\section*{Declaration of competing interest}

The authors declare that they have no known competing financial interests or personal relationships that could have
appeared to influence the work reported in this paper.

\section*{Acknowledgements} {The authors are grateful to the referees for
the helpful suggestions and comments. The Project is sponsored by ``Liaoning
BaiQianWan Talents Program". }

\vskip.2cm


\end{document}